\documentclass[12pt]{article}
\usepackage{amssymb}
\usepackage{latexsym}
\usepackage[all]{xy}
\def\mix{\mathbin{\diamondsuit}}
\newcounter{mytab}

\newtheorem{THM}{Theorem}[section]
\newtheorem{REM}{Remark}[section]

\newtheorem{LEM}{Lemma}[section]

\begin{document}
%
%%%%%%%%%%%%%%%%%%%%%%%%%%%%%%%%%%%%%%%%%%%%%%%%%%%%%%%%%%%%%%%%
%

\title{Constructions of Chiral Polytopes of Small Rank}

\author{Antonio Breda D'Azevedo\\
%\thanks{Supported in part by R\&D Unit {\it Mathem\'atica e Aplicacoes?}}\\
Departamento de Matem\'{a}tica\\
Universidade de Aveiro\\
P 3800 Aveiro, Portugal
\and and\\[.07in]
Gareth A. Jones\\
%\thanks{Supported in part ?}
School of Mathematics\\
University of Southampton\\
Southampton SO17 1BJ, United Kingdom
\and and\\[.07in]
Egon Schulte\thanks{Supported by NSA-grant H98230-07-1-0005 and NSF-grant  
DMS--0856675}\\
Department of Mathematics\\
Northeastern University\\
Boston, Massachusetts,  USA, 02115}

\date{ \today }
\maketitle

\begin{abstract}
An abstract polytope of rank $n$ is said to be {\em chiral\/} if its automorphism group has precisely two orbits on the flags, such that adjacent flags belong to distinct orbits. The present paper describes a general method for deriving new finite chiral polytopes from old finite chiral polytopes of the same rank. In particular, the technique is used to construct many new examples in ranks $3$, $4$ and $5$. 


\end{abstract}

\section{Introduction}
\label{intro}
Abstract polytopes are combinatorial structures with distinctive geometric, algebraic, or topological properties, in many ways more fascinating than convex polytopes and tessellations. The most symmetric structures, the abstract regular polytopes, generalize the traditional regular polytopes and regular tessellations and have been investigated extensively (see Coxeter~\cite{crp} and McMullen \& Schulte~\cite{arp}). The abstract regular polytopes of rank $3$ are (essentially) the maps on closed surfaces and have a long history of study across many areas of mathematics (see Coxeter \& Moser~\cite{cm} and Jones \& Singerman~\cite{js94}). 

By contrast, relatively little is known about (abstract) chiral polytopes, which form the most important class of nearly regular polytopes (see Schulte \& Weiss~\cite{SW1}). Chirality is a fascinating phenomenon which does not have a counterpart in the classical theory of traditional convex polytopes. Intuitively, a polytope of rank $n$ is chiral if it has maximal ``rotational symmetry" but lacks symmetry by ``reflection".  In rank $3$, the finite chiral polytopes are given by the irreflexible (chiral) maps on closed compact surfaces (see Coxeter \& Moser~\cite{cm}). The first family of chiral maps was constructed by Heffter~\cite{hef} in 1898 (see also Doro \& Wilson~\cite{doro}). It is well-known that there are infinitely many chiral maps on the $2$-torus (of genus 1), but for higher genus their appearance is rather sporadic, with the next example occurring on a surface of genus $7$ (see Wilson~\cite{wi}). On the other hand, for rank $n\geq4$ there are no chiral tessellations on compact euclidean space-forms of dimension $n-1$ (see Hartley, McMullen \& Schulte~\cite{hms}), so in particular, there are no chiral tessellations on the $(n-1)$-torus. (There do exist regular tessellations on the $(n-1)$-torus, called regular toroids, for any $n\geq 3$; see \cite[Sect. 6D,E]{arp}).

The quest for chiral polytopes of rank greater than $3$ has inspired much of the recent activity in this area. For a long time, finite examples were known only in rank $4$ (see \cite{Cox3,CW,MWei1990,nosc,SW2}). One particular type of construction in rank $4$ begins with a $3$-dimensional regular hyperbolic honeycomb and a faithful representation of its symmetry group as a group of complex M\"obius transformations (generated by the inversions in four circles cutting one another at the same angles as the corresponding reflection planes in hyperbolic space), and then derives chiral $4$-polytopes by applying modular reduction techniques to the corresponding matrix group (see Monson \& Schulte~\cite{monschmod} for a brief survey). 

For rank $5$, infinite chiral polytopes were found be applying the following general extension theorem to certain finite examples of rank $4$ (see \cite{SW3}). Any chiral polytope of rank $n$ with regular facets is itself the facet type of a chiral polytope of rank $n+1$; moreover, among all chiral polytopes of rank $n+1$ with facets isomorphic to the given chiral polytope of rank $n$, there exists a universal such polytope, whose automorphism group is a certain amalgamated product of the automorphism groups of the given polytope and its facet (see Section~\ref{ch5}).

Very recently, Conder, Hubard \& Pisanski~\cite{chp} succeeded in constructing the first known examples of finite chiral polytopes of rank $5$, by searching for normal subgroups of small index in the orientation preserving subgroups of certain $5$-generator Coxeter groups with string diagrams. 

In the present paper we describe a general method for deriving new chiral polytopes from old chiral polytopes of the same rank. This enables us to construct many new examples in ranks $3$, $4$ and $5$. The new polytopes have a group isomorphic to the direct product of the automorphism group of the old chiral polytope with the rotation (even) subgroup of the automorphism group of a certain regular polytope. The key idea behind this approach is an analogue of the mixing technique of \cite[Ch.7]{arp}.  

In Sections~\ref{polgr} and \ref{moreint} we review basic facts about regular and chiral polytopes and investigate the crucial intersection property of their automorphism groups. Section~\ref{chgroups} discusses chirality groups, an algebraic means of measuring the degree of irreflexibility of a chiral polytope.  
Section~\ref{monodchir} introduces mixing of groups and polytopes, and establishes criteria for when the mix of certain groups is a direct product. The four remaining sections then exploit the mixing technique to construct chiral polytopes in ranks $3$, $4$ and $5$.

\section{Polytopes and groups}
\label{polgr}

For general background material on abstract polytopes we refer the reader to \cite[Chs. 2,3]{arp}. Here we briefly review some basic concepts and terminology.

An (\emph{abstract\/}) \emph{polytope of rank\/} $n$, or simply an \emph{$n$-polytope\/}, is a partially ordered set $\mathcal{P}$ with a strictly monotone rank function with range $\{-1,0, \ldots, n\}$. An element of rank $j$ is a \emph{$j$-face\/} of $\mathcal{P}$, and a face of rank $0$, $1$ or $n-1$ is a \emph{vertex\/}, \emph{edge\/} or \emph{facet\/}, respectively. The maximal chains, or \emph{flags}, of $\mathcal{P}$ all contain exactly $n + 2$ faces, including a unique least face $F_{-1}$ (of rank $-1$) and a unique greatest face $F_n$ (of rank $n$). These faces $F_{-1}$ and $F_n$ are the \emph{improper} faces of $\mathcal{P}$; the other faces are the \emph{proper} faces of $\mathcal{P}$. Two flags are said to be \emph{adjacent} ($i$-\emph{adjacent}) if they differ in a single face (just their $i$-face, respectively). Then $\mathcal{P}$ is required to be \emph{strongly flag-connected}, in the sense that, if $\Phi$ and $\Psi$ are two flags, then they can be joined by a sequence of successively adjacent flags $\Phi = \Phi_0,\Phi_1,\ldots,\Phi_k = \Psi$, each containing $\Phi \cap \Psi$. Finally, $\mathcal{P}$ has the following homogeneity property (diamond condition):\ whenever $F \leq G$, with $F$ a $(j-1)$-face and $G$ a $(j+1)$-face for some $j$, then there are exactly two $j$-faces $H$ with $F \leq H \leq G$. This last property basically says that abstract polytopes are topologically ``real" and are close relatives of convex polytopes; in particular, this feature distinguishes abstract polytopes from other kinds of (thick) ranked incidence structures.

For any two faces $F$ of rank $j$ and $G$ of rank $k$ with $F \leq G$, we call
$G/F := \{ H \in \mathcal{P}\, | \, F \leq H \leq G \}$ a \emph{section} of $\mathcal{P}$; this is
a ($k-j-1$)-polytope in its own right. In particular, we can identify a face $F$ with the section $F/F_{-1}$. Moreover, $F_{n}/F$ is said to be the \emph{co-face at\/} $F$, or the \emph{vertex-figure at\/} $F$ if $F$ is a vertex.

We occasionally require more general ranked structures than polytopes. Following \cite[p.43]{arp}, a {\em pre-polytope\/} is a ranked partially ordered set sharing with polytopes all defining properties but one, namely strong flag-connectedness. Thus a pre-polytope is a polytope if and only if it is strongly flag-connected. A priori, a pre-polytope need not have any connectedness properties at all. However, in our applications, all pre-polytopes are at least {\em flag-connected\/}, meaning that, if $\Phi$ and $\Psi$ are two flags, they can be joined by a sequence of successively adjacent flags $\Phi = \Phi_0,\Phi_1,\ldots,\Phi_k = \Psi$, with no further assumptions on these flags. Thus strong flag-connectedness amounts to flag-connectedness of every section (including the pre-polytope itself). 

Now returning to polytopes, we say that a polytope $\mathcal{P}$ is \emph{regular\/} if its \textit{automorphism group\/} $\Gamma(\mathcal{P})$ (group of incidence preserving bijections) is transitive on the flags, and that $\mathcal{P}$ is \emph{chiral\/} if $\Gamma(\mathcal{P})$ has two flag orbits such that adjacent flags are always in distinct orbits. The group $\Gamma(\mathcal{P})$ of a regular or chiral polytope $\mathcal{P}$ has a well-behaved system of \emph{distinguished generators\/} obtained as follows.

If $\mathcal{P}$ is a regular $n$-polytope, then $\Gamma(\mathcal{P})$ is generated by involutions $\rho_0,\ldots,\rho_{n-1}$, where $\rho_i$ maps a fixed, or \emph{base\/}, flag $\Phi$ to the flag $\Phi^i$ $i$-adjacent to $\Phi$. These generators satisfy (at least) the standard Coxeter-type relations
\begin{equation}
\label{standardrel}
(\rho_i \rho_j)^{p_{ij}} = \epsilon \;\; \textrm{ for } i,j=0, \ldots,n-1,
\end{equation}
where $\epsilon$ denotes the identity element, and $p_{ii}=1$, $p_{ji} = p_{ij} =: p_{i+1}$ if $j=i+1$, and $p_{ij}=2$ otherwise. Note that the underlying Coxeter diagram is a string diagram. The numbers $p_j$ determine the (\emph{Schl\"afli}) 
\emph{type} $\{p_{1},\ldots,p_{n-1}\}$ of $\mathcal{P}$. Moreover, the group has the following {\em intersection property\/}:
\begin{equation}
\label{intprop}
\langle \rho_i \mid i \in I \rangle \cap \langle \rho_i \mid i \in J \rangle
= \langle \rho_i \mid i \in {I \cap J} \rangle
  \;\; \textrm{ for } I,J \subseteq \{0,1,\ldots,n-1\}.
\end{equation}
The elements
\[ \sigma_{i}:=\rho_{i-1}\rho_{i}\;\;\;(i=1,\ldots,n-1) \]
generate the {\em rotation subgroup\/} $\Gamma^{+}(\mathcal{P})$ of $\Gamma(\mathcal{P})$, which has index at most~$2$. We call $\mathcal{P}$ {\em directly regular\/} if this index is $2$.

A \emph{string C-group\/} is a group $\Gamma = \langle \rho_0,\ldots, \rho_{n-1} \rangle$ whose generators satisfy (\ref{standardrel}) and (\ref{intprop}); here, the ``C'' stands for ``Coxeter'', though not every C-group is a Coxeter group. The string C-groups are precisely the automorphism groups of regular polytopes, since, in a natural way, such a polytope can be uniquely reconstructed from $\Gamma$ (see \cite[\S 2E]{arp}). We often identify a regular polytope with its automorphism (string C-) group.

We write $[p_{1},p_{2},\ldots,p_{n-1}]$ for the Coxeter group whose underlying Coxeter diagram is a string with $n$ nodes and with $n-1$ branches labeled $p_{1},p_{2},\ldots,p_{n-1}$. This is the automorphism group of the \emph{universal} $n$-polytope $\{p_{1},\ldots,p_{n-1}\}$, which is indeed regular (see \cite[Sect. 3D]{arp}). 

All polytopes of rank $3$ (also called {\em polyhedra\/}) can be viewed as maps on surfaces, and all maps on surfaces satisfying the above homogeneity property are polytopes (of rank $3$). Note that a map on a surface has the homogeneity property if and only if its faces have more than one edge and do not self-touch (neither two vertices, nor two edges, of a face coincide). Recall from \cite{cm} that $\{p,q\}_r$ denotes the regular map obtained from the regular tessellation $\{p,q\}$ on the $2$-sphere or the euclidean or hyperbolic plane by identifying any two vertices $r$ steps apart along a Petrie polygon (a zigzag polygon along the edges such that any two, but no three, consecutive edges belong to a common face). Its group $[p,q]_r$ is the quotient of the Coxeter group $[p,q] = \langle\rho_{0},\rho_{1},\rho_{2}\rangle$ obtained by factoring out the single extra relation $(\rho_{0}\rho_{1}\rho_{2})^{r}=\epsilon$. In particular, $\{p,q\}_r$ is directly regular (or, equivalently, orientable) if and only if $r$ is an even integer. Note that here it is not implied that a polyhedron $\{p,q\}_r$ exists for all $p,q,r\geq 2$; for example, the group $[3,7]_{17}$ is trivial (see \cite[p.113]{cm}), so a polyhedron $\{3,7\}_{17}$ certainly does not exist.

If $\mathcal{P}$ is a chiral $n$-polytope, then $\Gamma(\mathcal{P})$ is generated by elements $\sigma_1,\ldots,\sigma_{n-1}$ associated with a base flag $\Phi=\{F_{-1},F_0,\ldots,F_{n}\}$ as follows. The generator $\sigma_i$ fixes the faces in $\Phi \setminus \{F_{i-1},F_{i}\}$ and cyclically permutes (``rotates") consecutive 
$i$-faces of $P$ in the (polygonal) section $F_{i+1}/F_{i-2}$ of rank $2$. By replacing a generator by its inverse if need be, we can further require that, if $F_{i}'$ denotes the $i$-face of $\mathcal{P}$ with $F_{i-1}<F_{i}'<F_{i+1}$ and $F_{i}'\neq F_i$, then $\sigma_{i}(F_{i}')=F_i$. The resulting generators $\sigma_1,\ldots,\sigma_{n-1}$ of $\Gamma(\mathcal{P})$ then satisfy (at least) the relations
\begin{equation}
\label{chiralrel}
\sigma_i^{p_i} = (\sigma_i\sigma_{i+1}\cdot\ldots\cdot\sigma_j)^{2} = \epsilon
  \;\; \textrm{ for } i,j=1,\dots,n-1,  \textrm{ with } i<j,
\end{equation}
where as before the numbers $p_i$ determine the {\it type\/} $\{p_1,\ldots,p_{n-1}\}$ of $\mathcal{P}$. The intersection property for the groups of chiral polytopes is more complicated than that for C-groups and can be described as follows (see \cite{SW1}). 

For $1\leq i\leq j\leq n-1$ define 
\[ \kappa_{i,j}:= \sigma_i\sigma_{i+1}\cdot\ldots\cdot\sigma_j, \]
and for $0\leq i\leq n$ let $\kappa_{0,i}:=\kappa_{i,n}:= \epsilon$; then $\kappa_{ii}=\sigma_{i}$ for $i\neq 0,n$. For $I\subseteq \{-1,0\ldots,n\}$ set
\[ \Gamma^{I} := \langle\kappa_{i,j}\mid i\leq j \mbox{ and } i-1,j\in I\rangle.  \]
Then the subgroup $\Gamma^{I}$ of $\Gamma(\mathcal{P})$ is trivial if $|I|\leq 1$; equals $\Gamma(\mathcal{P})$ if $I=\{-1,0,\ldots,n\}$; or is given by $\langle\sigma_{1},\ldots,\sigma_{n-2}\rangle$ or $\langle\{\sigma_{j}\mid j\neq i,i+1\}\cap\{\kappa_{i,i+1}\}\rangle$ if $I=\{-1,0,\ldots,n\}\setminus\{i\}$ with $i=n-1$ or $1\leq i\leq n-2$, respectively. Moreover, $\Gamma^{\{i-1,\ldots,j\}}=\langle\sigma_{i},\ldots,\sigma_{j}\rangle$ if $0\leq i\leq j\leq n$. The \emph{intersection property\/} for the group $\Gamma(\mathcal{P})$ of a chiral $n$-polytope $\mathcal{P}$ then takes the following form:
\begin{equation}
\label{intforchir}
\Gamma^{I} \cap \Gamma^{J} = \Gamma^{I\cap J}  
\;\; \textrm{ for } I,J \subseteq \{-1,0,\ldots,n\}.
\end{equation}
For polytopes of ranks $4$ and $5$, such as those dealt with in this paper, these conditions are equivalent to the following smaller sets of equalities of groups:\ for rank $4$, 
\begin{equation}
\label{chiralintprop4}
\langle\sigma_1\rangle \cap \langle\sigma_2\rangle = \langle\epsilon\rangle =
\langle\sigma_2\rangle \cap \langle\sigma_3\rangle, \quad
\langle\sigma_1,\sigma_2\rangle \cap \langle\sigma_2,\sigma_3\rangle
=\langle\sigma_2\rangle 
\end{equation}
(the intersection condition for rank $3$ amounts to the left-most equality in equation (\ref{chiralintprop4})); and for rank $5$, both sets (\ref{chiralintprop4}) and
\begin{equation}
\label{chiralintprop5}
\begin{array}{c}
\langle\sigma_1,\sigma_2,\sigma_3\rangle \cap 
\langle\sigma_2,\sigma_3,\sigma_4\rangle
= \langle\sigma_2,\sigma_3\rangle,\quad
\langle\sigma_1,\sigma_2,\sigma_3\rangle \cap \langle\sigma_3,\sigma_4\rangle
= \langle\sigma_3\rangle,\quad \\[.02in]
\langle\sigma_1,\sigma_2,\sigma_3\rangle \cap \langle\sigma_4\rangle
= \langle\epsilon\rangle .
\end{array}
\end{equation}
The equations in (\ref{chiralintprop4}) and (\ref{chiralintprop5}), respectively, correspond to the cases $n=4$ and $n=5$ of Lemma~\ref{fewint} described below. 

Note that the same intersection property (\ref{intforchir}) also holds for the rotation subgroup of a directly regular polytope. However, it might fail if the regular polytope is not directly regular. In fact, if $\mathcal{P}$ is a regular $4$-polytope whose facets and vertex-figures are not directly regular, then $\mathcal{P}$ is also not directly regular and
\[ \langle\sigma_{1},\sigma_{2}\rangle \cap \langle\sigma_{2},\sigma_{3}\rangle
= \langle\rho_{0},\rho_{1},\rho_{2}\rangle \cap \langle\rho_{1},\rho_{2},\rho_{3}\rangle
= \langle\rho_{1},\rho_{2}\rangle \,\gneqq\, \langle\sigma_{2}\rangle , \]
so the two subgoups on the left intersect in the dihedral group $ \langle\rho_{1},\rho_{2}\rangle$, not the cyclic group $\langle\sigma_{2}\rangle$. The $11$-cell $\{\{3,5\}_{5},\{5,3\}_{5}\}$ (see below for notation) discovered in \cite{chi,Gr2} is an example of a $4$-polytope of this kind. Its automorphism group is $L_{2}(11)$ and the facets and vertex-figures are hemi-icosahedra $\{3,5\}_{5}$ or hemi-dodecahedra $\{5,3\}_{5}$, respectively, both with group $A_5$; the two subgroups isomorphic to $A_5$ intersect in a dihedral group of order $10$.  (We use the ATLAS~\cite{atlas} notation $L_{2}(q)$ for the projective special linear group $PSL_{2}(q)$.)

Conversely, if $\Gamma$ is a group generated by elements $\sigma_1,\ldots,\sigma_{n-1}$ satisfying the relations (\ref{chiralrel}) and the respective intersection condition (\ref{intforchir}) in rank $n$ (that is, (\ref{chiralintprop4}) in rank $4$ and (\ref{chiralintprop5}) in rank $5$), then $\Gamma$ is the group of a chiral $n$-polytope, or the rotation subgroup (of index $2$) of the group of a directly regular polytope. In particular, the polytope is 
directly regular if and only if $\Gamma$ admits an involutory group automorphism mapping the set of generators $\sigma_1,\ldots,\sigma_{n-1}$ of $\Gamma$ to the new set of generators 
$\sigma_{1}^{-1},\sigma_{1}^{2}\sigma_2,\sigma_3,\ldots,\sigma_{n-1}$; 
in this case the automorphism is induced by conjugation with $\rho_0$ in $\Gamma(\mathcal{P})$ (and there are similar group automorphisms induced by conjugation with $\rho_i$ for any $i$). On the other hand, for a chiral polytope $\mathcal{P}$, the two flag orbits yield two sets of generators which are not conjugate in 
$\Gamma(\mathcal{P})$; thus a chiral polytope occurs in two \emph{enantiomorphic} (mirror image) forms.

If $\omega$ is any word in the generators $\sigma_1,\ldots,\sigma_{n-1}$ of $\Gamma$ (and their inverses), then the {\em enantiomorphic\/} (or {\em mirror image\/}) {\em word\/} $\overline{\omega}$ of $\omega$ is obtained from $\omega$ by replacing every occurrence of $\sigma_{1}$ by $\sigma_{1}^{-1}$ and $\sigma_2$ by $\sigma_{1}^{2}\sigma_2$, while keeping all $\sigma_j$ with $j\geq 3$ unchanged (and possibly applying some standard rules for inverses, powers, and cancellation). This definition is motivated by the characterization of direct regularity in the previous  paragraph. (Strictly speaking, its proper setting would be the free group on $n$ generators, but for simplicity we work here directly with $\Gamma$.)  For example, if $n=4$ and $\omega = \sigma_{2}\sigma_{3}^{-1}\sigma_{1}$, then $\overline{\omega}=\sigma_{1}^{2}\sigma_{2}\sigma_{3}^{-1}\sigma_{1}^{-1}$. If a word $\omega$ represents the trivial element in $\Gamma$, this may no longer be true for the new word $\overline{\omega}$ (unless $\Gamma$ admits a group automorphism as described). Note that $\overline{\omega}$ depends on the actual representation of an element of $\Gamma$ as a word $\omega$ in the generators, not only on the element itself. Moreover, observe that $\overline{\overline{\omega}} = \omega$ for all words $\omega$. 

We later exploit the following simple observation to establish the chirality of polytopes. An application is given in Section~\ref{ch5}.

\begin{LEM}
\label{wordcrit}
In the above situation, if $\Gamma$ is the rotation subgroup of a directly regular $n$-polytope and $\omega$ is any word in the distinguished generators of $\Gamma$, then the two elements of $\Gamma$ represented by $\omega$ and its enantiomorphic word $\overline{\omega}$ must necessarily have the same period. 
\end{LEM}

\noindent\textbf{Proof}.
The elements $\omega$ and $\overline{\omega}$ are conjugates by $\rho_0$ in the full automorphism group of the polytope, so they clearly have the same period.
\hfill $\square$

It is convenient to define $\Gamma^{+}(\mathcal{P}) := \Gamma(\mathcal{P})$ if $\mathcal{P}$ is a chiral polytope, so as to have common notation available for both chiral and directly regular polytopes.  Thus $\Gamma({\mathcal P})$, again considered with its distinguished generators, coincides with its rotation subgroup in this case. 

The simplest examples of chiral polytopes are the toroidal maps $\{4,4\}_{(b,c)}$, $\{3,6\}_{(b,c)}$ and $\{6,3\}_{(b,c)}$, with $b,c\neq 0$ and $b\neq c$ (see \cite{cm}). There seem to be fewer chiral polytopes than regular polytopes. For example, for chiral maps, the next occurrence (by genus) is on a surface of genus $7$; see \cite{maco} for a census of orientable chiral or regular maps of genus $2$ to $101$, and of non-orientable regular maps of genus $2$ to $202$. For rank $n\geq4$, there are no chiral tessellations on compact euclidean space-forms of dimension $n-1$ (see \cite{hms}), so the search of examples of chiral $n$-polytopes should begin with Schl\"afli symbols of hyperbolic type.

All the sections of a regular polytope are regular, and all the sections of a chiral polytope are either directly regular or chiral. In fact, for a chiral $n$-polytope, all the $(n-2)$-faces and all the co-faces at edges are actually directly regular; in particular, the $3$-faces of a chiral polytope of rank $5$ are directly regular.

If a regular or chiral $n$-polytope $\cal P$ has facets ${\cal P}_1$ and vertex-figures ${\cal P}_2$, we say that $\cal P$ is of {\em type\/} $\{{\cal P}_1,{\cal P}_2\}$ (this is a change of terminology from \cite{arp}).
By slight abuse of notation, the symbol $\{{\cal P}_1,{\cal P}_2\}$ will also denote the {\em universal\/} regular or chiral $n$-polytope of this type, provided any such polytopes exist at all (see \cite[p.97]{arp} and \cite[p.229]{SW2}). For example, the $11$-cell $\{\{3,5\}_{5},\{5,3\}_{5}\}$ is the universal regular $4$-polytope with hemi-icosahedral facets $\{3,5\}_{5}$ and hemi-dodecahedral vertex-figures $\{5,3\}_{5}$.

Recall that the \emph{order complex} of an $n$-polytope $\mathcal{P}$ is the (abstract) $(n-1)$-dimensional simplicial complex whose vertices are the proper faces of $\mathcal{P}$ (of ranks $0,\ldots,n-1$) and whose simplices are the chains (subsets of flags) which do not contain an improper face (of rank $-1$ or $n$). The order complex of a chiral polytope is orientable, and the order complex of a regular polytope is orientable if and only if the polytope is directly regular.

A polytope is {\em self-dual\/} if it admits a duality (incidence reversing bijection) onto itself.  A chiral polytope $\mathcal P$ is {\em properly self-dual\/} if it admits a duality preserving each of the two flag orbits under $\Gamma({\mathcal P})$ (that is, flags are mapped to flags in the same flag orbit); otherwise, $\mathcal P$ is {\em improperly self-dual\/} (see \cite{huwe,SW1}).

A self-dual regular, or properly self-dual chiral, $n$-polytope $\mathcal P$ always has a duality $\delta$ which fixes the base flag; this necessarily is a \textit{polarity} (duality of period $2$). If $\sigma_1,\ldots,\sigma_{n-1}$ are the distinguished generators of $\Gamma^{+}(\mathcal{P})$, then
\[ \delta\sigma_{j}\delta = \sigma_{n-j}^{-1} \qquad (j=1,\ldots,n-1) . \]
Thus conjugation by $\delta$ (in the group of all automorphisms and dualities of $\mathcal P$) induces an involutory group automorphism of $\Gamma^{+}(\mathcal{P})$. Conversely, if an $n$-polytope $\mathcal P$ is regular or chiral and $\Gamma^{+}(\mathcal{P})$ admits a group automorphism mapping $\sigma_{j}$ to 
$\sigma_{n-j}^{-1}$ for each $j$, then $\mathcal P$ is self-dual; in fact, $\mathcal{P}$ is properly self-dual if $\mathcal P$ is chiral. Note that an improperly self-dual chiral $n$-polytope need not possess a polarity, but it does if the rank $n$ is odd (see \cite{huwe}).

Let $\mathcal{P}$ and $\mathcal{Q}$ be two polytopes (or flag-connected pre-polytopes) of the same rank, not necessarily regular or chiral. Following \cite[p.43]{arp}, a mapping $\gamma: \mathcal{P}\rightarrow\mathcal{Q}$ is called a {\em covering\/} if it preserves incidence of faces, ranks of faces, and adjacency of flags; then $\gamma$ is necessarily surjective, by the flag-connectedness of $\mathcal{Q}$. We say that $\mathcal{P}$ {\em covers\/} $\mathcal{Q}$ if there exists a covering $\gamma: \mathcal{P}\rightarrow\mathcal{Q}$.
 
A common way to obtain coverings of polytopes is by the construction of quotients. Let $\mathcal{P}$ be an $n$-polytope, let $N$ be a subgroup of $\Gamma(\mathcal{P})$, and let $\mathcal{P}/N$ denote the set of $N$-orbits on $\mathcal{P}$. If $F$ is a face of $\mathcal{P}$, we write $N\cdot F$ for its orbit under $N$. On $\mathcal{P}/N$ we can introduce a partial order as follows:\ if $\widehat{F},\widehat{G}\in \mathcal{P}/N$, then $\widehat{F}\leq \widehat{G}$ if and only if $\widehat{F} = N\cdot F$ and $\widehat{G} = N\cdot G$ for some faces $F$ and $G$ of $\mathcal{P}$ with $F\leq G$. The set $\mathcal{P}/N$ together with this partial ordering is called the {\em quotient of $\mathcal{P}$ with respect to $N$}.  A quotient of a polytope need not be a polytope (in fact, it need not even be a pre-polytope).

The chiral polytopes in this paper are obtained from other chiral polytopes by {\em mixing\/} them with a directly regular polytope. To begin with, let $\Delta=\langle\alpha_{1},\ldots,\alpha_{k}\rangle$ and 
$\Delta'=\langle\alpha_{1}',\ldots,\alpha_{k}'\rangle$ be two groups, each with $k$ specified generators (here we could allow a generator to be trivial, although we will not encounter this case). If we define
\[  \beta_{j} := (\alpha_{j},\alpha_{j}') \in
\Delta\times \Delta'\,  \quad (j = 1,\ldots,k), \]
then the subgroup
\begin{equation}
\label{mixdef}
\Delta\mix\Delta' :=
\langle\beta_{1},\ldots,\beta_{k}\rangle
\end{equation}
of $\Delta\times \Delta'$ is called the {\em mix\/} of (the {\em component groups\/}) $\Delta$ and $\Delta'$ (see \cite[Ch.7A]{arp}). In our applications, $\Delta$ and $\Delta'$ will be the rotation subgroups of chiral or directly regular polytopes. 

More specifically, let $\mathcal{P}$ and $\mathcal{Q}$ be chiral or regular $n$-polytopes with (rotation sub-) groups 
$\Gamma^{+}(\mathcal{P}) =  \langle\sigma_{1},\ldots,\sigma_{n-1}\rangle$ and
$\Gamma^{+}(\mathcal{Q}) = \langle\sigma_{1}',\ldots,\sigma_{n-1}'\rangle$,
respectively.  (Here we represent a regular polytope by its rotation subgroup rather than the full automorphism group. We do not explicitly require that a regular polytope be directly regular, but in our applications this will always be the case.)  Now the mix of $\Gamma^{+}(\mathcal{P})$ and $\Gamma^{+}(\mathcal{Q})$ is given by
\begin{equation}
\label{mixgr}
\Gamma^{+}(\mathcal{P}) \mix \Gamma^{+}(\mathcal{Q}) :=
\langle\tau_{1},\ldots,\tau_{n-1}\rangle,
\end{equation}
where
\[  \tau_{j} := (\sigma_{j},\sigma_{j}')  \quad
(j = 1,\ldots,n-1). \]
This certainly is a group satisfying relations as in (\ref{chiralrel}), but in general it will not have the intersection property (\ref{intforchir}) with respect to its generators $\tau_{1},\ldots,\tau_{n-1}$. If, however, the intersection property is satisfied, then $\Gamma^{+}(\mathcal{P}) \mix \Gamma^{+}(\mathcal{Q})$ is either the group of a chiral $n$-polytope or the rotation subgroup of a directly regular $n$-polytope, denoted $\mathcal{P}\mix \mathcal{Q}$ and called the {\em mix} of $\mathcal{P}$ and $\mathcal{Q}$.

Note that the mix of two polytopes is an analogue of the parallel product (or join) construction occurring in the study of maps and hypermaps (see \cite{ant1,ant2,orb,wilson}). 

\section{More on the intersection property}
\label{moreint}

In practice, the verification of the intersection property (\ref{intforchir}) for a group $\Gamma$ can often be reduced to the consideration of only a few intersections of subgroups, as described in the following lemma given without proof (see \cite[p.511]{SW1}). The cases $n=4$ and $n=5$ were illustrated earlier.

\begin{LEM}
\label{fewint}
Let $n\geq4$, and let $\Gamma :=\langle\sigma_1,\ldots,\sigma_{n-1}\rangle$ be a group satisfying (\ref{chiralrel}). Suppose that the subgroup $\Gamma_{n-1}:=\langle\sigma_1,\ldots,\sigma_{n-2}\rangle$ of 
$\Gamma$ has the intersection property (\ref{intforchir}), and that the following intersection conditions hold:
\begin{equation}
\label{inductint}
\Gamma_{n-1} \cap \langle\sigma_j,\ldots,\sigma_{n-1}\rangle
= \langle\sigma_j,\ldots,\sigma_{n-2}\rangle
\qquad (j=2,\ldots,n-1) .
\end{equation}
Then $\Gamma$ itself has the intersection property (\ref{intforchir}).
\end{LEM}

Next we establish, in Lemma~\ref{quotcrit}, an analogue for chiral polytopes of the quotient criterion of \cite[p.56]{arp} for regular polytopes. We refer to it again as the {\em quotient criterion\/}. This often enables us to determine if a given group with presentation (\ref{chiralrel}) has the intersection property (\ref{intforchir}). 

\begin{LEM}
\label{quotcrit}
Let $\Gamma :=\langle\sigma_1,\ldots,\sigma_{n-1}\rangle$ be a group satisfying (\ref{chiralrel}), and let $\Lambda :=\langle\lambda_1,\ldots,\lambda_{n-1}\rangle$ be a group satisfying (\ref{chiralrel}) and the intersection property (\ref{intforchir})\  (or equivalently, (\ref{chiralintprop4}) or (\ref{chiralintprop5}) when $n=4$ or $5$, respectively). If the mapping $\sigma_j \rightarrow \lambda_j$ for $j=1,\ldots,n-1$ induces a homomorphism $\pi: \Gamma \rightarrow \Lambda$, which is one-to-one on 
$\Gamma_{n-1}:=\langle\sigma_1,\ldots,\sigma_{n-2}\rangle$ or on 
$\Gamma_{0}:=\langle\sigma_2,\ldots,\sigma_{n-1}\rangle$, 
then $\Gamma$ also has the intersection property (and $\pi$ induces a covering between the $n$-polytopes associated with $\Gamma$ and $\Lambda$).
\end{LEM}

\noindent\textbf{Proof}.
The proof is similar to the proof of \cite[Thm. 2E17]{arp}. Suppose that $\pi$ is one-to-one on $\Gamma_{n-1}$, so that, in particular, $\Gamma_{n-1}$ has the intersection property. By the previous lemma it suffices to verify that
\[ \Gamma_{n-1} \cap \langle\sigma_j,\ldots,\sigma_{n-1}\rangle =
\langle\sigma_j,\ldots,\sigma_{n-2}\rangle
\quad (j=2,\ldots,n-1) . \]
Let $\varphi \in \Gamma_{n-1} \cap \langle\sigma_j,\ldots,\sigma_{n-1}\rangle$. Then
\[ \pi(\varphi) \in
\langle\lambda_1,\ldots,\lambda_{n-2}\rangle \cap
\langle\lambda_j,\ldots,\lambda_{n-1}\rangle =
\langle\lambda_j,\ldots,\lambda_{n-2}\rangle ,\]
since $\Lambda$ has the intersection property. Hence $\pi(\varphi)$ has a pre-image in
$\langle\sigma_j,\ldots,\sigma_{n-2}\rangle$. This must necessarily coincide with $\varphi$, since $\pi$ is one-to-one on $\Gamma_{n-1}$. Thus $\varphi\in\langle\sigma_j,\ldots,\sigma_{n-2}\rangle$. It follows 
that $\Gamma$ has the intersection property.
\hfill $\square$

Concluding this section we deal with the particularly interesting case when the quotient criterion guarantees the ``polytopality" of the mix $\Gamma^{+}(\mathcal{P}) \mix \Gamma^{+}(\mathcal{Q})$ for chiral or directly regular polytopes $\mathcal P$ and $\mathcal Q$. Clearly, given $\mathcal P$ and $\mathcal Q$ the canonical projections
\[ \pi_{P}: \Gamma^{+}(\mathcal{P})\,\mix\, \Gamma^{+}(\mathcal{Q}) 
\rightarrow \Gamma^{+}(\mathcal{P})\;\, \mbox{ and }\;\,
\pi_{Q}: \Gamma^{+}(\mathcal{P})\,\mix\, \Gamma^{+}(\mathcal{Q}) 
\rightarrow \Gamma^{+}(\mathcal{Q}) \]
of the mix onto its component groups are surjective homomorphisms. If the vertex-figures (or facets, respectively) of $\mathcal P$ and $\mathcal Q$ are isomorphic, then the restrictions of $\pi_P$ and $\pi_Q$ to the subgroup $\langle\tau_2,\ldots,\tau_{n-1}\rangle$ (or $\langle\tau_1,\ldots,\tau_{n-2}\rangle$, respectively) are one-to-one. Hence Lemma~\ref{quotcrit} applies and proves that the mix has the 
intersection property. This establishes the following result.

\begin{LEM}
\label{isomvert}
Let $\mathcal{P}$ and $\mathcal{Q}$ be chiral or directly regular $n$-polytopes with isomorphic vertex-figures (or facets, respectively). Then $\Gamma^{+}(\mathcal{P})\mix\Gamma^{+}(\mathcal{Q})$ is the group of a chiral or directly regular $n$-polytope, whose vertex-figures (or facets, respectively) are isomorphic to those of $\mathcal P$ and $\mathcal Q$.
\end{LEM}

The conditions on $\mathcal{P}$ and $\mathcal{Q}$ described in the next two sections will ensure that the resulting polytope is indeed chiral.

\begin{REM}
\label{remint}
{\rm(a)} Inspection of the argument shows that, for the mix $\Gamma^{+}(\mathcal{P})\mix\Gamma^{+}(\mathcal{Q})$ in Lemma~\ref{isomvert} to have the intersection property (\ref{intforchir}), it suffices to assume that only one of its component groups has this property, but still requiring isomorphic facets or vertex-figures, respectively, for the two structures. In particular, it is sufficient to require that only one component, $\mathcal{P}$ or $\mathcal{Q}$, is a chiral or directly regular $n$-polytope while the other component still is a flag-connected pre-polytope of rank $n$. (The groups of chiral or directly regular pre-polytopes do not in general have the intersection property.)  In other words, in the above, both $\mathcal P$ and $\mathcal Q$ must be flag-connected pre-polytopes of rank $n$ but only one of them actually needs to be a polytope of the appropriate kind. \\
{\rm(b)} Generalizing in yet another direction, it also is enough to assume that the vertex-figures (or facets, respectively) of one of the polytopes $\mathcal P$ and $\mathcal Q$ cover those of the other (they need not actually be isomorphic, though in our applications this is typically the case).
\end{REM}

\section{Chirality groups}
\label{chgroups}

In this section we find it convenient to represent $n$-polytopes as quotients of the universal $n$-polytope $\mathcal{U} = \{\infty,\ldots,\infty\}$ associated with the (universal) string Coxeter group $W=W_n$ on $n$ generators $r_{0},\ldots,r_{n-1}$. Thus $W=W_n$ has the presentation
\begin{equation}
\label{univrel}
r_{i}^{2} = (r_{i}r_{j})^{2} = \epsilon \;\; \textrm{ for } i,j
=0,\dots,n-1, \textrm{ with } i<j-1.
\end{equation}
This is a string C-group and its rotation subgroup $W^{+}=W_{n}^{+}$ of index $2$ is generated by
$s_{i}:=r_{i-1}r_{i}\;(i=1,\ldots,n-1)$, with a presentation given by
\begin{equation}
\label{univrotrel}
(s_{i}s_{i+1} \ldots s_{j})^{2} = \epsilon
  \;\; \textrm{ for } i,j=1,\dots,n-1,  \textrm{ with } i<j.
\end{equation}

We now describe the concept of the chirality group of an abstract polytope, an invariant which in some sense measures the degree of mirror asymmetry (irreflexibility) of the polytope; this is an analogue of the chirality group for hypermaps introduced in \cite{bjns} (see also \cite{bbn}). Here we restrict ourselves to chiral or directly regular polytopes, although the concept can be generalized to more general classes of polytopes. Our discussion is in terms of automorphism groups of polytopes, not monodromy groups as in \cite{bjns}, but the two approaches are equivalent at least for chiral or directly regular polytopes (in general, for other kinds of polytopes, the definitions would have to involve the monodromy group --- see \cite{ha1,huor}). For ease of presentation we restrict ourselves to polytopes, although the concept of chirality group applies more generally to pre-polytopes (see Remark~\ref{remint}).

Let $\mathcal{P}$ be a chiral or directly regular $n$-polytope, let $\Gamma^{+}(\mathcal{P})$ be its (rotation) group, and let $M$ be the normal subgroup of $W^{+}$ such that $\Gamma^{+}(\mathcal{P}) = W^{+}/M$. Here, $\mathcal{P}$ is directly regular if and only if $M$ is also normal in $W$, that is, if and only if $M^{r_{0}}=M$. (Throughout we write $M^{g}$ for $gMg^{-1}$.)  Given the polytope $\mathcal P$ the largest normal subgroup of $W$ contained in $M$ is
\[ M_{W}:= M \cap M^{r_{0}} , \]
and the smallest normal subgroup of $W$ containing $M$ is
\[ M^{W}:= MM^{r_{0}}. \]
Both are also normal subgroups of $W^+$. This relationship is displayed in the following diagram:
\begin{equation}
\label{display}
\begin{array}{c}
\begin{picture}(180,85)
\put(71,70){$M^W$}
\put(67,65){\line(-1,-1){15}}
\put(88,65){\line(1,-1){15}}
\put(45,36){$M$}
\put(98,36){$M^{r_{0}}$}
\put(71,0){$M_W$}
\put(67,16){\line(-1,1){15}}
\put(88,16){\line(1,1){15}}
\end{picture}
\end{array}
\end{equation}
The corresponding quotients $\mathcal{P}_{W}:=\mathcal{U}/M_W$ and $\mathcal{P}^{W}:=\mathcal{U}/M^W$ of $\mathcal U$ are ``polytope-like" ranked partially ordered sets which may or may not be polytopes. (In fact, when $M^{W}=W^{+}$ the quotient $\mathcal{P}^W$ consists only of two flags, so certainly $\mathcal{P}^W$ is not even a pre-polytope.)  Informally we can think of $\mathcal{P}_{W}$ as the smallest ``regular" (reflexible) cover of $\mathcal{P}$, and of $\mathcal{P}^{W}$ as the largest (polytope-like) ``regular" ranked partially ordered set covered by $\mathcal{P}$; once again, with no claim of polytopality. In any case, if $\mathcal{P}$ is finite, then $\mathcal{P}_{W}$ and $\mathcal{P}^{W}$ are also finite, since both $M_W$ and $M^W$ have finite index in $W$. Recall here that the intersection of two subgroups of finite index also has finite index.

The following result is a direct analogue of \cite[Prop. 2]{bjns}.

\begin{LEM}
\label{fourgroups}
Any two of the four groups $M^{W}/M$, $M/M_{W}$, $M^{W}/M^{r_{0}}$ and $M^{r_{0}}/M_{W}$ are isomorphic.
\end{LEM}

\noindent\textbf{Proof}.  We have
\[ M^{W}/M = MM^{r_{0}}/M \simeq M^{r_{0}}/(M\cap M^{r_{0}}) =
M^{r_{0}}/M_{W},\]
and similarly $M^{W}/M^{r_{0}} \simeq M/M_{W}$. Moreover, conjugation by $r_0$ induces the isomorphisms $M^{W}/M \simeq M^{W}/M^{r_{0}}$ and $M/M_{W} \simeq M^{r_{0}}/M_{W}$.
\hfill $\square$

We will call the common quotient group of Lemma~\ref{fourgroups} the {\em chirality group\/} $X(\mathcal{P})$ of $\mathcal{P}$, and its order the {\em chirality index\/} $\kappa = \kappa(\mathcal{P})$ of $\mathcal{P}$. Note that $X(\mathcal{P})$ is trivial if and only if $\mathcal{P}$ is directly regular. Thus $X(\mathcal{P})$ measures algebraically how far the polytope $\mathcal{P}$ deviates from being reflexible. Note that, since $X(\mathcal{P})=M^{W}/M$, we may view $X(\mathcal{P})$ as a normal subgroup of $\Gamma^{+} (\mathcal{P}) = W^{+}/M$ with quotient 
\begin{equation}
\label{sugbringamma}
\Gamma^{+}(\mathcal{P})/X(\mathcal{P}) = (W^{+}/M)/(M^{W}/M) = W^{+}/M^{W}.
\end{equation}
The most extreme case of chirality (irreflexibility) occurs when $X(\mathcal{P})$ and $\Gamma^{+}(\mathcal{P})$ are the same; that is, when $M^{W}=W^+$. In this case the polytope $\mathcal{P}$ is said to be {\em totally chiral\/}. Every chiral polytope $\mathcal{P}$ whose automorphism group is simple is an example of a totally chiral polytope; in fact, $M^W/M$, being a non-trivial normal subgroup of $W^{+}/M$, must coincide with the full group, giving $M^{W}=W^+$. However, not every totally chiral polytope has a simple group as automorphism group (see \cite{bjns} for examples of rank $3$). In Section~\ref{ch5} we also meet an example of a chiral $5$-polytope whose group is $S_6$ (and hence almost simple), but which is not totally chiral.

The next lemma says that, if $\mathcal{P}$ is totally chiral, then the rotation subgroup of the smallest regular cover $\mathcal{P}_{W}$ of $\mathcal{P}$ is the direct product of the automorphism group of $\mathcal{P}$ with itself.

\begin{LEM}
\label{ant}
In the above notation, if $\mathcal{P}$ is totally chiral, then
\[ W^{+}/M_{W} \;\cong\; \Gamma^{+}(\mathcal{P}) \times  \Gamma^{+}(\mathcal{P}) 
=   \Gamma(\mathcal{P}) \times  \Gamma(\mathcal{P}). \]
\end{LEM}

\noindent\textbf{Proof}.
Since $\mathcal{P}$ is totally chiral, we have $W^{+}=M^{W}=MM^{r_0}$ and hence
\[ W^{+}/M_{W} \,=\, MM^{r_0}/M_{W} \,=\, M/M_{W}\cdot M^{r_0}/M_{W}.\]
The two factors on the right have trivial intersection, since $M\cap M^{r_0} = M_{W}$; hence the internal product is direct. Moreover,
\[ M/M_{W} = M/(M\cap M^{r_0}) \cong MM^{r_0}/M^{r_0} = W^{+}/M^{r_0} \cong W^{+}/M 
= \Gamma^{+}(\mathcal{P}) = \Gamma (\mathcal{P}), \]
and similarly, $M^{r_0}/M_{W} \cong \Gamma(\mathcal{P})$. Now the lemma follows.
\hfill $\square$
\medskip

It is usually very difficult to compute the chirality group $X(\mathcal{P})$ of a chiral polytope $\mathcal{P}$ by hand (however, see \cite{bbn} for the explicit computation of the chirality groups of the chiral torus maps). We conclude this section with some general remarks about generating sets of $X(\mathcal{P})$ which sometimes allow the computation of $X(\mathcal{P})$ with GAP~\cite{gap}. In particular, this enabled us to calculate the chirality group of a 5-polytope in Section~\ref{ch5}. 

Consider the canonical epimorphism 
\[ \Gamma(\mathcal{P}) = W^{+}/M \longrightarrow 
\Gamma(\mathcal{P})/X(\mathcal{P}) = (W^{+}/M)/(M^{W}/M) = W^{+}/M^{W} ,\]
whose kernel is $M^{W}/M = X(\mathcal{P})$ (see (\ref{sugbringamma})). Since $W^+$ is generated by $s_{1},\ldots,s_{n-1}$, we can write 
\[ \Gamma(\mathcal{P}) = \langle s_{1},\ldots,s_{n-1} \mid \mathcal{R}\rangle ,\]
where $\mathcal{R}$ is the set of relators in the generators $s_{1},\ldots,s_{n-1}$ that defines $\Gamma(\mathcal{P})$ as a quotient of $W^+$ by $M$; in other words, $M$ is the normal closure of $\mathcal{R}$ in $W^+$, denoted $\langle\mathcal{R}\rangle^{W^{+}}$. Since $W^{+}$ is normal in $W$, we have
\[ M^{r_{0}} = (\langle\mathcal{R}\rangle^{W^{+}})^{r_0} = \langle\mathcal{R}^{r_0}\rangle^{W^{+}},\]
where $\mathcal{R}^{r_0}$ denotes the set of relators $R^{r_0}$ (say) in the generators $s_{1},\ldots,s_{n-1}$ obtained from the relators $R$ in $\mathcal{R}$ by replacing each generator $s_{j}$ in $R$ by the product 
\[ (s_{1}s_{2}\cdots s_{j-1})(s_{1}s_{2}\cdots s_{j})^{-1}.\] 
This latter product is $s_{1}^{-1}$, $s_{1}^{2}s_2$, or $s_j$, according as $j=1$, $j=2$, or $j\geq 3$.
Now, considering the image of $M^{r_{0}} = \langle\mathcal{R}^{r_0}\rangle^{W^{+}}$ under the canonical projection $W^{+}\rightarrow W^{+}/M$, we see that the chirality group $X(\mathcal{P}) = MM^{r_0}/M$ is the subgroup of $\Gamma(\mathcal{P})=W^{+}/M$ given by 
\[ X(\mathcal{P})=\langle\mathcal{R}^{r_0}\rangle^{\Gamma(\mathcal{P})} ,\]
the normal closure in $\Gamma(\mathcal{P})$ of the image of $\mathcal{R}^{r_0}$ in $\Gamma(\mathcal{P})$, which for simplicity we have again denoted by $\mathcal{R}^{r_0}$. Thus $X(\mathcal{P})$ is the normal closure of $\mathcal{R}^{r_0}$ in $\Gamma(\mathcal{P})$. 

If $\mathcal{P}$ is a finite chiral polytope and the defining relators $\mathcal{R}$ for its group $\Gamma(\mathcal{P})$ are explicitly known, the calculation of the corresponding normal closure in $\Gamma(\mathcal{P})$ often becomes manageable in GAP~\cite{gap} and then outputs the chirality group of the polytope.

\section{Mixing}
\label{monodchir}

We now describe the basic construction that lies at the heart of our discovery of new chiral polytopes. Suppose that $M$ and $K$, respectively,  are the normal subgroups of $W^{+}=\langle s_{1},\ldots,s_{n-1}\rangle$ associated with two $n$-polytopes $\mathcal{P}$ and $\mathcal{Q}$, where $\mathcal{P}$ is chiral and $\mathcal{Q}$ is directly regular; this latter condition says that $K$ is normal in $W$, but $M$ is not. Our methods apply more generally in situations where the quotients $W^{+}/M$ and $W^{+}/K$ may not necessarily be associated with polytopes (just with pre-polytopes --- see Remark~\ref{remone}), but here, for the sake of simplicity, we will make this assumption. The situation is depicted in the following diagram:
\begin{equation}
\label{dia}
\begin{array}{c}
\begin{picture}(180,145)
\put(70,132){$W^{+}$}
\put(67,125){\line(-1,-1){35}}
\put(88,125){\line(1,-1){35}}
\put(77.5,90){\line(0,1){35}}
\put(18,77){$M$}
\put(124,77){$M^{r_{0}}$}
\put(61,77){$K\!=\!K^{r_0}$}
\put(67,68){\line(-1,-1){35}}
\put(88,68){\line(1,-1){35}}
\put(22,33){\line(0,1){35}}
\put(132,33){\line(0,1){35}}
\put(-40,19){$M\cap K \,=\, N$}
\put(124,19){$N^{r_{0}} \,=\,(M\cap K)^{r_{0}}$}
\end{picture}
\end{array}
\end{equation}
Here $N:=M\cap K$ is a normal subgroup of $W^+$ which may or may not be normal in $W$. We later describe conditions guaranteeing that $N$ is not normal in $W$.

Writing $\sigma_{i} := s_{i}M$ and $\sigma_{i}' := s_{i}K$ for $i=1,\ldots,n-1$, we then have
\[ \begin{array}{lclcl}
W^{+}/M = \langle \sigma_{1},\ldots,\sigma_{n-1} \rangle =
\Gamma^{+}(\mathcal{P}), \\[.05in]
W^{+}/K = \langle \sigma_{1}',\ldots,\sigma_{n-1}' \rangle =
\Gamma^{+}(\mathcal{Q}) ,
\end{array} \]
where $\sigma_{1},\ldots,\sigma_{n-1}$ and $\sigma_{1}',\ldots,\sigma_{n-1}'$ are the two sets of distinguished generators. We now investigate the mix
\begin{equation}
\label{mixmono}
W^{+}/M \,\mix\, W^{+}/K  \;=\; \langle \tau_{1},\ldots,\tau_{n-1}\rangle ,
\end{equation}
the subgroup of the direct product $W^{+}/M \,\times\, W^{+}/K$ generated by $\tau_{i}:=(\sigma_{i},\sigma_{i}')$ for $i=1,\ldots,n-1$.

First observe that the natural projections onto the two components define surjective homomorphisms of the mix in (\ref{mixmono}) onto the component groups $W^{+}/M$ and $W^{+}/K$, so the latter are quotients of the mix. Moreover, since both sets of generators $\sigma_i$ of $W^{+}/M$ and $\sigma_{i}'$ of $W^{+}/K$ satisfy the defining relations (\ref{univrotrel}) for $W^+$, these relations continue to hold for the generators $\tau_{i}$ of $W^{+}/M \,\mix\, W^{+}/K$, so the latter is necessarily a quotient group of $W^{+}$, under a surjective homomorphism that maps $w\in W^{+}$ to $(wM,wK)$ and hence has kernel $M\cap K =N$. Thus,
\begin{equation}
\label{mixmk}
W^{+}/N \,\simeq\, W^{+}/M \,\mix\, W^{+}/K ,
\end{equation}
under the isomorphism given by $wN \rightarrow (wM,wK)$. We now have the following lemma.

\begin{LEM}
\label{mixintersect}
In the above situation of diagram~(\ref{dia}), if the mix $W^{+}/M \,\mix\, W^{+}/K$ has the intersection property with respect to its generators $\tau_1,\ldots,\tau_{n-1}$, then it is the rotation group of a chiral or directly regular $n$-polytope.
\end{LEM}

Next we make additional assumptions on the subgroups $M$ and $K$ of $W^+$ which allow us to explicitly determine the mix; in fact, the mix turns out be the full direct product.

\begin{LEM}
\label{mixdir}
In the above situation of diagram~(\ref{dia}), suppose that the groups $W^{+}/M$ and $W^{+}/K$ have no nontrivial common quotient. (For example, this holds if one of the groups $W^{+}/M$ or $W^{+}/K$ is simple and the other does not have it as a quotient.)
Then,\\
\noindent {\rm (a)} $W^{+} = MK$,\\
\noindent {\rm (b)} $W^{+}/N \,\simeq\, M/N \times K/N$,\\
\noindent {\rm (c)} $W^{+}/M \,\mix\, W^{+}/K \;=\; W^{+}/M
\,\times\, W^{+}/K$.
\end{LEM}

\noindent\textbf{Proof}.
For part (a) note that $MK$ is normal in $W^{+}$ and
\[ (W^{+}/M)/(MK/M)\, \simeq\, W^{+}/MK \,\simeq\, (W^{+}/K)/(MK/K) .\]
Hence $W^{+}/MK$ is a nontrivial common quotient of $W^{+}/M$ and $W^{+}/K$, except when $MK=W^{+}$. This establishes part (a). Then part (b) follows from
\[ W^{+}/N = MK/N = (M/N)\cdot (K/N)\]
and the observation that $(M/N) \cap (K/N)$ is trivial since $N=M\cap K$. Finally, from part~(a) we obtain
\[ W^{+}/M = MK/M \simeq M/(M\cap K) = M/N, \]
and similarly $W^{+}/K \simeq  K/N$, so that now part (c) is implied by part~(b) and equation (\ref{mixmk}).
\hfill $\square$

We mention in passing that the strong assumption of Lemma~\ref{mixdir} that there are no nontrivial common quotients, is not satisfied in general if $\mathcal{P}$ is a totally chiral polytope with a non-simple automorphism group $W^{+}/M$.

Restated in terms of polytopes, our previous discussion yields the following lemma.

\begin{LEM}
\label{polmix}
Let $\mathcal{P}$ and $\mathcal{Q}$ be $n$-polytopes, let $\mathcal{P}$ be chiral, and let $\mathcal{Q}$ be directly regular. Suppose that the groups $\Gamma^{+}(\mathcal{P}) \;(= \Gamma(\mathcal{P}))$ and $\Gamma^{+}(\mathcal{Q})$ have no nontrivial common quotient. (For example, this holds if one of the groups $\Gamma^{+}(\mathcal{P})$ or $\Gamma^{+}(\mathcal{Q})$ is simple and the other does not have it as a quotient.) Then,
\[  \Gamma^{+}(\mathcal{P})\,\mix\, \Gamma^{+}(\mathcal{Q}) \;=\;
\Gamma^{+}(\mathcal{P})\,\times\, \Gamma^{+}(\mathcal{Q}) .\]
Moreover, if this group has the intersection property with respect to its generators $\tau_1,\ldots,\tau_{n-1}$, then it is the rotation group of a chiral or directly regular $n$-polytope, the mix $\mathcal{P}\mix\mathcal{Q}$ of $\mathcal{P}$ and $\mathcal{Q}$.
\end{LEM}

In the interesting special case when the vertex-figures (or facets) of the two component polytopes are isomorphic we can appeal to Lemma~\ref{isomvert} to immediately arrive at the following lemma.

\begin{LEM}
\label{isofacvert}
Let $\mathcal{P}$ and $\mathcal{Q}$ be $n$-polytopes with isomorphic vertex-figures (or facets, respectively), let $\mathcal{P}$ be chiral, and let $\mathcal{Q}$ be directly regular. Suppose that the groups $\Gamma^{+}(\mathcal{P})$ and $\Gamma^{+}(\mathcal{Q})$ have no nontrivial common quotient. Then $\Gamma^{+}(\mathcal{P})\times\Gamma^{+}(\mathcal{Q})$ is the rotation group of a chiral or directly regular $n$-polytope, the mix $\mathcal{P}\mix \mathcal{Q}$. Moreover, the vertex-figures (or facets, respectively) are isomorphic to those of $\mathcal P$ and $\mathcal Q$.
\end{LEM}

It remains to settle the critical question of chirality of the mix $\mathcal{P}\mix \mathcal{Q}$. Clearly, one possibility is to invoke the criterion of Lemma~\ref{wordcrit} for enantiomorphic pairs of words. However, as an alternative, we can also establish a criterion based on chirality groups. This can be done as follows. As before, let $\mathcal{P}$ be chiral, and let $\mathcal{Q}$ be directly regular. Clearly, $\mathcal{P}\mix \mathcal{Q}$ is not chiral in general; in fact, if $\mathcal{Q}$ covers $\mathcal{P}$, then $\mathcal{P}\mix \mathcal{Q}$ is isomorphic to $\mathcal{Q}$ and hence is regular. Ignoring the question of polytopality for a moment, for $\mathcal{P}\mix \mathcal{Q}$ to be chiral we need to prove that $X(\mathcal{P}\mix \mathcal{Q})$ is not the trivial group.

\begin{LEM}
\label{detchgr}
If $\mathcal{P}$ is a chiral $n$-polytope and $\mathcal{Q}$ is a directly regular $n$-polytope, then $X(\mathcal{P}\mix \mathcal{Q})$ is a normal subgroup of $X(\mathcal{P})$ and 
\[ X(\mathcal{P}\mix \mathcal{Q}) \cong (M^{r_0}\cap K)M/M. \]
Moreover, $X(\mathcal{P}\mix \mathcal{Q})$ is the trivial group if and only if $M\cap K \leq M^{r_0}$. In particular, if $\mathcal{P}\mix \mathcal{Q}$ is directly regular, then $\mathcal{P}\mix \mathcal{Q}$ covers the smallest regular cover $\mathcal{P}_W$ of $\mathcal{P}$.
\end{LEM}

\noindent\textbf{Proof}.
Recall from (\ref{mixmk}) that $\mathcal{P}\mix \mathcal{Q}$ is associated with the normal subgroup $N=M\cap K$ of $W^+$. Then,
\[ X(\mathcal{P}\mix \mathcal{Q}) = N^{r_0}/N_W
= (M^{r_0} \cap K)/ (M^{r_0} \cap K \cap M), \]
bearing in mind that $K^{r_0}=K$. Now consider the homomorphism
\begin{equation}
\label{homgamma}
\begin{array}{lrlll}
\gamma: & M^{r_0} \cap K & \longrightarrow & M^{r_0}M/M  = X(\mathcal{P}) \\
               & w & \longrightarrow & w M.
\end{array}
\end{equation}
Then the kernel and image of $\gamma$ are given by $Ker(\gamma) = M^{r_0} \cap K \cap M$ and
\[ Im(\gamma) = \{wM\mid w\in M^{r_0} \cap K \} = (M^{r_0} \cap K)M/M . \]
Hence
\[ X(\mathcal{P}\mix \mathcal{Q})  = (M^{r_0} \cap K)/ (M^{r_0} \cap K \cap M)
= (M^{r_0} \cap K)/Ker(\gamma) \cong Im(\gamma) = (M^{r_0} \cap K)M/M, \]
where the group on the right is a normal subgroup of $M^{r_0}M/M = X(\mathcal{P})$. This proves the first part of the lemma.

For the second part note that $(M^{r_0}\cap K)M/M$ is trivial if and only if $M^{r_0}\cap K \leq  M$, or equivalently, if and only if $M\cap K \leq M^{r_0}$. In other words, $X(\mathcal{P}\mix \mathcal{Q})$ is trivial if and only if the normal subgroup $N = M\cap K$ associated with $\mathcal{P}\mix \mathcal{Q}$ lies in the normal subgroup $M_W$ associated with the smallest regular cover $\mathcal{P}_W$ of $\mathcal{P}$. However, the latter implies that $\mathcal{P}\mix \mathcal{Q}$ covers $\mathcal{P}_W$.
\hfill $\square$

\begin{REM}
\label{remnext}
More generally, if $\mathcal{P}$ and $\mathcal{Q}$ are arbitrary chiral or directly regular $n$-polytopes (with no condition on one being directly regular), the chirality group of the mix $\mathcal{P}\mix \mathcal{Q}$ embeds into the direct product of the chirality groups of the components; that is,
\[ X(\mathcal{P}\mix \mathcal{Q}) \leq X(\mathcal{P}) \times X(\mathcal{Q}) .\]
For the proof of this more general result, replace $\gamma$ in (\ref{homgamma}) by the homomorphism
\[\begin{array}{lrlll}
\gamma: & M^{r_0} \cap K^{r_0} & \longrightarrow & M^{r_0}M/M \times K^{r_0}K/K
= X(\mathcal{P}) \times X(\mathcal{Q})  \\
               & w & \longrightarrow & (wM, wK),
\end{array} \]
but otherwise follow the same line of argument. When $\mathcal{P}$ is chiral and $\mathcal{Q}$ is directly regular (and hence $X(\mathcal{Q})$ is trivial), we recover Lemma~\ref{detchgr}.
\end{REM}

Next we employ Lemmas~\ref{ant} and \ref{detchgr} to establish the chirality of the mix $\mathcal{P}\mix \mathcal{Q}$ under the assumptions that $\mathcal{P}$ is (chiral and) totally chiral and $\mathcal{Q}$ is directly regular. Now suppose for a moment that the corresponding mix $\mathcal{P}\mix \mathcal{Q}$ is directly regular in this case. Then $X(\mathcal{P}\mix \mathcal{Q})$ is trivial, and hence $\mathcal{P}\mix \mathcal{Q}$ covers $\mathcal{P}_{W}$ by Lemma~\ref{detchgr}. On the other hand, by Lemma~\ref{ant}, since $\mathcal{P}$ is totally chiral, the group $W^{+}/M_W$ of $\mathcal{P}_W$  is isomorphic to $\Gamma^{+}(\mathcal{P})\times \Gamma^{+}(\mathcal{P})$. It follows that $\Gamma^{+}(\mathcal{P})\times \Gamma^{+}(\mathcal{P})$ must be a quotient of the group $\Gamma^{+}(\mathcal{P})\mix \Gamma^{+}(\mathcal{Q})$ of $\mathcal{P}\mix \mathcal{Q}$. Thus, if both polytopes $\mathcal{P}$ and $\mathcal{Q}$ are finite, then
\[ |\Gamma^{+}(\mathcal{P})\mix \Gamma^{+}(\mathcal{Q})| \equiv 0\; \bmod |\Gamma^{+}(\mathcal{P})|^2 . \]
On the other hand, $\Gamma^{+}(\mathcal{P})\mix \Gamma^{+}(\mathcal{Q})$ is a subgroup of $\Gamma^{+}(\mathcal{P})\times \Gamma^{+}(\mathcal{Q})$ and so its order divides $|\Gamma^{+}(\mathcal{P})|\cdot|\Gamma^{+}(\mathcal{Q})|$. Hence, $ |\Gamma^{+}(\mathcal{P})|\cdot|\Gamma^{+}(\mathcal{Q})| \equiv 0\,\bmod |\Gamma^{+}(\mathcal{P})|^2 $, 
or equivalently, 
\[ |\Gamma^{+}(\mathcal{Q})|\equiv 0\; \bmod |\Gamma^{+}(\mathcal{P})| .\] 
This proves the following lemma.

\begin{LEM}
\label{chircritone}
Let $\mathcal{P}$ be chiral and totally chiral, let $\mathcal{Q}$ be directly regular, and suppose that $\Gamma^{+}(\mathcal{P})\mix\Gamma^{+}(\mathcal{Q})$ has the intersection property. If
$\Gamma^{+}(\mathcal{P})\times\Gamma^{+}(\mathcal{Q})$ does not have a subgroup which has a quotient $\Gamma^{+}(\mathcal{P})\times\Gamma^{+}(\mathcal{P})$ (in particular if $\mathcal{P}$ and $\mathcal{Q}$ are finite and $|\Gamma^{+}(\mathcal{P})|$ does not divide $|\Gamma^{+}(\mathcal{Q})|$), then $\mathcal{P}\mix\mathcal{Q}$ is a chiral polytope.
\end{LEM}

\begin{REM}
\label{remone}
As indicated earlier, some of our results carry over to certain normal subgroups $M$ and $K$ of $W^+$ whose quotient groups $W^{+}/M$ and $W^{+}/K$ are not necessarily associated with polytopes. Such quotients always are associated with pre-polytopes (see Remark~\ref{remint}). Still, there are interesting examples where the mix itself is polytopal even though one of its components is not. We will meet examples in later sections (see Theorems~\ref{locspher535} and~\ref{locspher353}).
\end{REM}

\section{Chiral polyhedra}
\label{polyh}

Polyhedra (polytopes of rank $3$) and maps on surfaces have been studied for well over 100 years, and deep connections with other branches of mathematics have been discovered, including hyperbolic geometry, Riemann surfaces, number fields and Galois theory (see \cite{cm,js94}). In the past few years there has been great progress in the computer-aided enumeration of regular or chiral maps by genus, leading to the creation of a complete census of regular or chiral maps on orientable surfaces of genus up to 101 (see \cite{maco}); this approach is based on new algorithms for finding low index normal  subgroups. Moreover, recently the set of possible genera (known as the genus spectrum) for regular or chiral maps has attracted a lot of attention (for example, see \cite{breda,masi}). 

In this section we briefly illustrate our method by mixing two polyhedra of type $\{3,7\}$ to construct infinite sequences of chiral polyhedra of type $\{3,7\}$. Recall that a {\em Hurwitz group\/} is a group generated by two elements of orders $2$ and $3$ whose product has order $7$ (yielding $(2,3,7)$ generating triples). The Hurwitz groups are precisely the rotation subgroups of chiral or directly regular polyhedra of type $\{3,7\}$ (see \cite{hur} and \cite[p. 264]{js87}). (For a given polyhedron, we may take $\sigma_{1}\sigma_2$ and $\sigma_{1}^{-1}$ as generators; the appropriate intersection condition holds since $3$ and $7$ are coprime.) It is known that the projective special linear groups $L_{2}(q)$, $q$ a prime power, are Hurwitz groups precisely when $q=7$, or when $q=p$ for any prime $p \equiv \pm 1 \bmod 7$, or when $q=p^3$ for any prime $p\not \equiv 0,\pm 1 \bmod 7$, but for no other values of $q$ (see \cite{hur,mac}); by Dirichlet's Theorem on primes in arithmetic progressions, this gives an infinite list of values of $q$. In fact, by a theorem of Macbeath~\cite{mac}, the group $L_2(q)$ yields three Hurwitz surfaces, and three polyhedra, corresponding to three orbits of $Aut(L_{2}(q))$ on $(2,3,7)$ generating triples, if $q = p \equiv \pm 1 \bmod 7$, and one if $q = 7$ or $q = p^3$ for $p\not\equiv 0,\pm 1 \bmod 7$. Moreover, the corresponding polyhedra $\mathcal{Q}$ (say) of type $\{3,7\}$ cannot be chiral and hence must necessarily be directly regular; in fact, if $L_{2}(q)$ is the quotient of $W^{+}$ by a normal subgroup $N$, then $N$ must necessarily be normal in $W$ (see \cite[Cor. 9]{bjns}). 

For the chiral component of the mix we take a polyhedron $\mathcal{P}$ of type $\{3,7\}$ whose group is a Ree group $Re(3^{f})= {}^{2}G_{2}(3^f)$; the latter is defined for all odd $f\geq 1$, has order $3^{3f}(3^{3f}+1)(3^{f}-1)$, and is simple when $f>1$. The Ree groups were first described in \cite{Ree1,Ree2}, and a good account of their properties can also be found in \cite[Ch. XI]{hup}. It was pointed out in the proof of \cite[Thm. 12]{bjns} that $Re(3^{f})$ with $f>1$ is a Hurwitz group yielding a totally chiral hypermap; this hypermap is in fact a chiral polyhedron of type $\{3,7\}$ which is totally chiral. Now, mixing $\mathcal{P}$ with a directly regular polyhedron $\mathcal{Q}$ arising from a group $L_{2}(q)$ by considering it as a Hurwitz group (possibly in several non-equivalent ways), gives a chiral polyhedron of type $\{3,7\}$ with group $Re(3^{f})\times L_{2}(q)$, as long as $f>1$ (see Lemmas~\ref{isofacvert} and \ref{chircritone}). 
In fact, the crucial condition of Lemma~\ref{chircritone} is satisfied:\ the group $\Gamma^{+}(\mathcal{P})\times\Gamma^{+}(\mathcal{Q}) = Re(3^{f})\times L_{2}(q)$ does not have a subgroup with quotient $\Gamma^{+}(\mathcal{P})\times\Gamma^{+}(\mathcal{P}) = Re(3^{f})\times Re(3^{f})$. Otherwise, we could 
project on to the second direct factor to obtain a subgroup of $L_2(q)$ with the Ree group as quotient; however, the subgroups of the groups $L_2(q)$ are known for all $q$ and the only non-abelian composition factors which occur are isomorphic to $A_5$ or $L_2(q')$ where $q$ is a power of $q'$ (see \cite[Ch. XII]{Dic} and \cite[Ch. II]{huppert}). Thus, the direct products $Re(3^{f})\times L_{2}(q)$ are Hurwitz groups associated with chiral polyhedra. For example, when $q=7$, any odd $f>1$ gives a chiral polyhedron of type $\{3,7\}$ with automorphism group $Re(3^{f})\times L_{2}(7)$.

A similar construction is possible using the Suzuki groups $Sz(2^f) = {}^2B_2(2^f)$, introduced by Suzuki in \cite{Suz} (see also \cite[Ch.~XI]{hup}). These are simple groups of order $2^{2f}(2^{2f}+1)(2^f-1)$ for odd $f>1$. It is shown in \cite[\S 6]{JS} that each Suzuki group is the automorphism group of a chiral map of type $\{4,5\}$, which can be regarded as a totally chiral polyhedron $\mathcal P$ of this type. As in the case of the Hurwitz groups, the techniques developed by Macbeath in \cite{mac} show that there are infinitely many values of $q$ for which $L_2(q)$ is the automorphism group of a directly regular polyhedron $\mathcal Q$ of type $\{4,5\}$: for instance one can take $q$ to be any prime $p \equiv \pm 1$ or $\pm 9$ mod~$40$ (again, by Dirichlet's Theorem, there are infinitely many such primes). Our construction then yields a chiral polyhedron of type $\{4,5\}$ with automorphism group $Sz(2^f) \times L_2(q)$ for each odd $f>1$.

\section{Chiral 4-polytopes with spherical vertex-figures}
\label{sph}

In this section we employ the mixing technique to construct finite chiral $4$-polytopes with spherical or toroidal facets and with spherical vertex-figures. This yields examples of locally spherical or locally toroidal $4$-polytopes. Recall that a $4$-polytope is said to be {\em locally spherical\/} if all its facets and vertex-figures are spherical maps, and {\em locally toroidal\/} if all its facets and vertex-figures are spherical or toroidal maps, with at least one toroidal. The regular or chiral toroidal maps are all of the forms $\{4,4\}_{(b,c)}$, $\{6,3\}_{(b,c)}$ or $\{3,6\}_{(b,c)}$ (see \cite{cm}); we write $[4,4]^{+}_{(b,c)}$, $[6,3]^{+}_{(b,c)}$ or $[3,6]^{+}_{(b,c)}$ for the corresponding rotation subgroup, which is the full automorphism group if the map is chiral.

As input we take a finite chiral $4$-polytope $\mathcal P$ and a finite regular $4$-polytope $\mathcal Q$, both of the same Schl\"afli type $\{r,s,t\}$ and both already locally spherical or locally toroidal. Since there are no chiral $4$-polytopes of (global) spherical or euclidean type, $\{r,s,t\}$ must necessarily be a hyperbolic Schl\"afli symbol (see \cite{hms} and \cite[Sect. 6H]{arp}). The admissible Schl\"afli symbols (with spherical vertex-figure) are
\[ \{r,s,t\} \,=\, \{3,5,3\}, \{5,3,4\}, \{5,3,5\}, \{4,4,3\}, 
\{6,3,3\},\{6,3,4\},\{6,3,5\} . \]
The vertex-figures of both polytopes $\cal P$ and $\cal Q$ are Platonic solids $\{s,t\}$. Hence Lemma~\ref{isofacvert} applies, provided the corresponding groups $\Gamma^{+}({\cal P})$ and $\Gamma^{+}({\cal Q})$ have no common non-trivial quotient. The latter condition is guaranteed if one of these groups is simple and the other does not have it as a quotient. In our applications, the simple group will usually come from the chiral polytope; then this polytope is also totally chiral, so the mix $\mathcal{P}\mix \mathcal{Q}$ is chiral by Lemma~\ref{chircritone} if $|\Gamma^{+}(\mathcal{P})|$ does not divide $|\Gamma^{+}(\mathcal{Q})|$.

We describe some particularly interesting examples but do not attempt to fully exhaust all possibilities. We mainly concentrate on examples where the simple group is a group $L_2(p)$ with $p$ an odd prime. There are many variants of our constructions. 

Recall that $L_{2}(p)$ has order $p(p^{2}-1)/2$ if $p$ is an odd prime. In our theorems we will usually make the assumption that, for a directly regular polytope $\mathcal{Q}$, the order $|\Gamma({\mathcal Q})|$ of its full automorphism group is not divisible by $p(p^{2}-1)$, or equivalently, that $|\Gamma^{+}({\mathcal Q})|$ is not divisible by $p(p^{2}-1)/2$. Note that this implies in particular that $L_{2}(p)$ cannot be a quotient of $\Gamma^{+}({\mathcal Q})$. This non-divisibility condition enables us to appeal directly to Lemma~\ref{chircritone}, although the lemma frequently applies in more general situations as well. 

\smallskip
\noindent
{\bf Locally spherical polytopes\/}

We begin with locally spherical chiral polytopes.

\begin{THM}
\label{locspher535}
Let $p$ be a prime, let $p\equiv\pm 1\bmod 5$ or $p=5$, and let $\mathcal Q$ be a finite locally spherical directly regular $4$-polytope of type $\{5,3,5\}$ such that $p(p^{2}-1)$ does not divide $|\Gamma(\mathcal{Q})|$. Then there exists a locally spherical chiral $4$-polytope of type $\{5,3,5\}$ with group $L_{2}(p)\times\Gamma^{+}({\mathcal Q})$.
\end{THM}

\noindent\textbf{Proof}.
We know from \cite[Thm. A]{JoLo} that under our assumptions on $p$ there are normal subgroups of $[5,3,5]^+$, not normal in $[5,3,5]$, whose quotient is isomorphic to $L_2(p)$. The corresponding quotients may not all be associated with polytopes, but if they are, then the corresponding polytopes $\mathcal P$ are chiral. Nevertheless, as indicated in Remark~\ref{remone}, even though one component group may not be polytopal, we still are able to apply our method if the intersection property can be verified for the resulting mix. In the present case, the mix 
\[ \Gamma:= L_{2}(p)\times \Gamma^{+}({\cal Q})=\langle\tau_1,\tau_2,\tau_3\rangle \] 
necessarily has its subgroups $\Gamma_{3}:=\langle\tau_1,\tau_2\rangle$ and 
$\Gamma_{1}:=\langle\tau_2,\tau_3\rangle$ isomorphic to $[5,3]^+$ and $[3,5]^+$, respectively,
since the latter are the rotation subgroups of the facet and vertex-figure groups for $\mathcal Q$ and certainly have the corresponding subgroups of $L_2(p)$ as homomorphic images. Thus the quotient criterion applies (via the projection onto the component group $\Gamma^{+}({\cal Q})$) and establishes that the mix has the intersection property, yielding a chiral polytope. Moreover, the mix is the direct product $L_{2}(p)\times\Gamma^{+}({\mathcal P})$, since the component groups have no nontrivial common quotient groups.
\hfill $\square$

For instance, let $\mathcal Q$ be the classical regular star-polytope $\{\frac{5}{2},3,5\}$ (of type $\{5,3,5\}$) in euclidean $4$-space (see \cite{crp}); here the fractional entry in the Schl\"afli symbol indicates that its $2$-faces are pentagrams $\{\frac{5}{2}\}$. Its combinatorial automorphism group $\Gamma(\mathcal{Q})$ is isomorphic to its geometric symmetry group given by the Coxeter group $H_{4} = [5,3,3]$ of order $14\,400$ (with suitably chosen distinguished generators), and $\Gamma^{+}(\mathcal{Q})=H_{4}^+$. Moreover, $\Gamma(\mathcal{Q})$ is obtained from $[5,3,5]$ by imposing the single extra relation 
\[ (\rho_{0}\rho_{1}\rho_{2}\rho_{3}\rho_{2}\rho_{1})^{3}=\epsilon \]
(see \cite{mgrsp} or \cite[p. 213]{arp}). Hence Theorem~\ref{locspher535} applies (if $p\neq 5)$ and yields locally spherical chiral $4$-polytopes of type $\{5,3,5\}$ with groups isomorphic to $L_{2}(p)\times H_{4}^+$. Note here that $H_{4}^{+}$ is a central product of two copies of $SL_{2}(5)$, so it has a central subgroup $C_2$ with quotient $L_{2}(5)^{2} \cong A_{5}^{2}$. There are many other possible choices for $\mathcal Q$ (for example, see \cite{halee,monschmod}); however, we cannot choose Coxeter's $57$-cell $\{\{5,3\}_5,\{5,3\}_5\}$, since this is not directly regular (see \cite{chd}).

\begin{THM}
\label{locspher353}
Let $p$ be a prime with $p\equiv\pm 1\bmod 5$. Let $p\equiv 1,3,4,5 \mbox{ or }\, 9\bmod 11$ and $3\pm 2\sqrt{5}$ both be squares $\bmod\;p$, or let $p\equiv 2,6,7,8 \mbox{ or }\, 10\bmod 11$. Let $\mathcal Q$ be a finite locally spherical directly regular $4$-polytope of type $\{3,5,3\}$ such that $p(p^{2}-1)$ does not divide $|\Gamma(\mathcal{Q})|$. Then there exists a locally spherical chiral $4$-polytope of type $\{3,5,3\}$ with group $L_{2}(p)\times\Gamma^{+}({\mathcal Q})$.
\end{THM}

\noindent\textbf{Proof}.
We now appeal to the results in \cite[Thm. 1.3]{JoLoMed}. In particular, under our assumptions on $p$ there exist normal subgroups of $[3,5,3]^+$, not normal in $[3,5,3]$, whose quotient is isomorphic to $L_2(p)$. Hence we can proceed as in the proof of Theorem~\ref{locspher535}, but now with $[3,5,3]$ in place of $[5,3,5]$ and with $\Gamma_{3}= [3,5]^+$ and $\Gamma_{1}= [5,3]^+$; recall that $[3,5]^{+}\cong [5,3]^{+}\cong A_5$. Thus the mix is again a chiral polytope with the desired properties. 
\hfill $\square$

As star-polytopes of type $\{3,5,3\}$ do not exist (see \cite{crp}) and the $11$-cell $\{\{3,5\}_5,\{5,3\}_5\}$ is not directly regular, the search for appropriate candidates for $\mathcal Q$ in Theorem~\ref{locspher353} is slightly more involved. Interesting examples of directly regular polytopes of type $\{3,5,3\}$ (and $\{5,3,5\}$) were constructed in \cite{monschmod} by applying modular reduction techniques, with moduli given by primes in $\mathbb{Z}[\tau]$ (with $\tau$ the golden ratio); their automorphism groups are given by certain finite orthogonal groups. Any of these polytopes can serve for $\mathcal Q$. For example, if $q$ is an odd (rational) prime with $q \equiv \pm 2 \bmod{5}$, then the automorphism group of the polytope is $O_1(4,q^2,\epsilon)$ with
\[ \epsilon :=\left\{ 
\begin{array}{cl}
+1, & \mbox{ if }  q \equiv 3, 12, 23, 27, 37, 
38, 42, 47, 48, 53 \bmod{55};\\
-1, & \mbox{ if } q \equiv 2, 7, 8, 13, 17, 18, 28, 32, 43, 52 \bmod{55}.
\end{array}
\right. \]
Recall here that, if $O(4,q^2,\epsilon)$ denotes the full orthogonal group of a $4$-dimensional non-singular orthogonal vector space over $GF(q^2)$ (with Witt index $2$ or $1$ according as $\epsilon=+1$ or $-1$), then $O_1(4,q^2,\epsilon)$ is the subgroup of index $2$ in $O(4,q^2,\epsilon)$ generated by the reflections whose spinor norm is $1$ (see \cite{art,monsch}). When applied with these polytopes as directly regular components of the mix, Theorem~\ref{locspher353} gives locally spherical chiral $4$-polytopes of type $\{3,5,3\}$ with groups 
\[ L_{2}(p)\times O_{1}^{+}(4,q^2,\epsilon) \cong L_{2}(p)\times \Omega(4,q^2,\epsilon),\] 
at least when $p\!\nmid\! (q^{2}\pm 1)$; here, $O_{1}^{+}(4,q^2,\epsilon)$ stands for the even subgroup of $O_{1}(4,q^2,\epsilon)$, which is isomorphic to the commutator subgroup $\Omega(4,q^2,\epsilon)$ of 
$O(4,q^2,\epsilon)$ (for this isomorphism, see \cite[Thms. 5,14, 5.17]{art} and \cite[p.297]{monsch}, and observe that the elements of $O_{1}^{+}(4,q^2,\epsilon)$ all have determinant $1$ and spinor norm $1$). Note that $O_{1}(4,q^2,\epsilon)$ has order $q^{2}(q^{2}-1)(q^{2}-\epsilon)$, so $p(p^{2}-1)$ certainly does not divide $|\Gamma(\mathcal{Q})|$ when $p\!\nmid\! (q^{2}\pm 1)$ (our conditions on $p$ and $q$ imply that $p\neq q$). 

\begin{THM}
\label{locspher534}
Let $p$ be a prime with $p \equiv \pm 1\mbox{ or }\pm 9 \bmod 40$, with either $1+\sqrt 5$ or $1-\sqrt 5$ a square $\bmod\; p$. Let $\mathcal Q$ be a finite locally spherical directly regular $4$-polytope of type $\{5,3,4\}$ such that $p(p^{2}-1)$ does not divide $|\Gamma(\mathcal{Q})|$. Then there exists a locally spherical chiral $4$-polytope of type $\{5,3,4\}$ with group $L_{2}(p)\times\Gamma^{+}({\mathcal Q})$.
\end{THM}

\noindent\textbf{Proof}.
As in the two previous cases, our assumptions on $p$ guarantee the existence of normal subgroups $M$ of $[5,3,4]^+$, not normal in $[5,3,4]$, whose quotient is $L_2(p)$. Here we arrived at the conditions on $p$ by carefully inspecting the calculations for \cite[Thm. 4.5(3)]{Lo}, and imitating the argument for $[3,5,3]$ in \cite{JoLoMed}. The congruence condition on $p$ (needed to give elements of order $4$ and $5$ in $L_2(p)$) implies that $p \equiv \pm 1 \bmod 5$, so quadratic reciprocity implies that $5$ is a square $\bmod\; p$. There seems to be no simple answer for when an element $1\pm\sqrt 5$ is a square $\bmod\; p$. However, $(1+\sqrt 5)(1-\sqrt 5)=-4$, so if $p\equiv 1\bmod 4$ then $-1$ is a square $\bmod\; p$ and hence either both or neither of $1\pm\sqrt 5$ are squares, while if $p\equiv -1 \bmod 4$ then $-1$ is not a square $\bmod\; p$, so exactly one of $1\pm\sqrt 5$ is a square $\bmod\; p$. (For example, if $p=31$ then $\sqrt 5 = \pm 6$, so $1\pm\sqrt 5 = 7$ or $-5$, and quadratic reciprocity shows that $7$ is a square $\bmod\; 31$, whereas $-5$ is not.)  Thus we always get such subgroups $M$ if $p \equiv -1 \mbox{ or } -9 \bmod 40$. (Note that similar remarks also apply to the elements $3\pm 2\sqrt 5$ occurring in Theorem~\ref{locspher353}.)

The construction of the corresponding chiral polytopes then proceeds as before, now with $\Gamma_{3}= [5,3]^+$ and $\Gamma_{1}= [3,4]^{+}\cong S_4$. 
\hfill $\square$

Interesting directly regular $4$-polytopes $\mathcal Q$ of type $\{5,3,4\}$ with automorphism groups $D_{s}^{6}\rtimes [3,5]$, $s\geq 2$, are obtained from \cite[Thm. 8C5]{arp}; here, $D_s$ denotes the dihedral group of order $2s$. More precisely, they arise as the duals of certain polytopes $2^{\mathcal{K},\mathcal{G}(s)}$, namely when $\mathcal{K}$ is the icosahedron $\{3,5\}$ and $\mathcal{G}(s)$ is the Coxeter diagram whose nodes are the vertices of $\{3,5\}$, and whose branches connect pairs of antipodal vertices of $\{3,5\}$ (but no other pairs of vertices) and are marked $s$. Thus, when $p \equiv -1 \mbox{ or } -9 \bmod 40$ and $p\nmid s$, Theorem~\ref{locspher534} and its proof provide locally spherical chiral $4$-polytopes of type $\{5,3,4\}$ with groups $L_{2}(p)\times (D_{s}^{6}\rtimes [3,5])^+$ (with an appropriate interpretation of even subgroup).

\smallskip
\noindent
{\bf Locally toroidal polytopes\/}

Next we turn to locally toroidal chiral $4$-polytopes. Interesting examples with groups $L_{2}(p)$ as automorphism group were constructed in \cite{SW2} and we may choose them here for the chiral component $\mathcal P$. The congruence classes of primes $p$ that actually occur depend on the underlying Schl\"{a}fli type; this explains the conditions on $p$ occurring in our theorems. Any directly regular $4$-polytope of the specified kind can serve as components $\mathcal Q$; finite examples are readily available for each Schl\"afli type (for example, see \cite[Chs. 10, 11]{arp} and \cite{mshf,mstor,moschtwo,monwei2,SW2}). Again we do not strive for the most general results. We begin with the type $\{4,4,3\}$. 

\begin{THM}
\label{loctor443}
Let $p$ be a prime with $p \equiv 1 \bmod 8$, let $b$ and $c$ be positive integers such that $p=b^{2}+c^{2}$, let $m\geq 2$ be an integer, and let $\mathcal Q$ be a finite directly regular $4$-polytope of type $\{\{4,4\}_{(m,0)},\{4,3\}\}$ such that $p(p^{2}-1)$ does not divide $|\Gamma(\mathcal{Q})|$. Then there exists a chiral $4$-polytope of type $\{\{4,4\}_{(m,0)},\{4,3\}\}$ if $p\!\mid\!m$ or $\{\{4,4\}_{(mb,mc)},\{4,3\}\}$ if $p\!\nmid\! m$, whose group is $L_{2}(p)\times \Gamma^{+}({\cal Q})$. The facets themselves are regular or chiral, respectively, in the two cases.
\end{THM}

\noindent\textbf{Proof}.
We know from \cite[Cor. 7.7]{SW2} that, under our assumptions on $p$, $b$ and $c$, there exists a chiral $4$-polytope $\mathcal P$ of type $\{\{4,4\}_{(b,c)},\{4,3\}\}$ with group $L_2(p)$. Now Lemmas~\ref{isofacvert} and \ref{chircritone} apply and show that $\mathcal{P}\mix\mathcal{Q}$ is a chiral $4$-polytope with vertex-figures $\{4,3\}$ and group $L_{2}(p)\times \Gamma^{+}({\cal Q})$. To determine the structure of the facets, again let $\tau_{i}=(\sigma_{i},\sigma_{i}')$ for $i=1,2,3$ denote the generators of 
\[ \Gamma^{+}({\cal P}) \mix\Gamma^{+}({\cal Q}) = 
L_{2}(p)\times \Gamma^{+}({\cal Q}),\]
where $\sigma_1,\sigma_2,\sigma_3$ generate $\Gamma^{+}({\cal P})$ and $\sigma_{1}',\sigma_{2}',\sigma_{3}'$ generate $\Gamma^{+}({\cal Q})$. We know that the facets of $\mathcal{P}\mix\mathcal{Q}$ are isomorphic to the mix $\{4,4\}_{(b,c)} \mix \{4,4\}_{(m,0)}$. Now recall that $[4,4]^{+}_{(b,c)} \cong C_p\rtimes C_4$ and $[4,4]^{+}_{(m,0)} \cong (C_m\times C_m)\rtimes C_4$, where in each case the first factor of the semi-direct product is given by the translation subgroup of the map and the second factor $C_4$ by the $4$-fold rotation about the base vertex. Here two cases can occur. 

First, if $p\!\mid\!m$, then $\{4,4\}_{(m,0)}$ covers $\{4,4\}_{(b,c)}$, so the facets are isomorphic to $\{4,4\}_{(m,0)}$. Hence $\mathcal{P}\mix \mathcal{Q}$ is of type $\{\{4,4\}_{(m,0)},\{4,3\}\}$ in this case. Second, suppose that $p\!\nmid\! m$, so $p$ and $m$ are relatively prime. Then $[4,4]^{+}_{(b,c)} \mix [4,4]^{+}_{(m,0)}$ is a subgroup of $[4,4]^{+}_{(b,c)} \times [4,4]^{+}_{(m,0)}$ isomorphic to $(C_p \times C_m\times C_m)\rtimes C_4$. But the latter is precisely the group $[4,4]^{+}_{(mb,mc)}$ of order $4pm^{2}$. Thus the facets are isomorphic to $\{4,4\}_{(mb,mc)}$ in this case and $\mathcal{P}\mix \mathcal{Q}$ is 
of type $\{\{4,4\}_{(mb,mc)},\{4,3\}\}$.
\hfill $\square$

As an example, if $\mathcal{Q}$ is the universal $4$-polytope $\{\{4,4\}_{(3,0)},\{4,3\}\}$ with full automorphism group $S_6\times C_2$ (see \cite[p.366]{arp}), then Theorem~\ref{loctor443} yields chiral $4$-polytopes of type $\{\{4,4\}_{(3b,3c)},\{4,3\}\}$ with groups $L_{2}(p)\times A_{6}\times C_{2}$.

\begin{THM}
\label{loctor633}
Let $p$ be a prime with $p \equiv 1 \bmod 12$, let $b$ and $c$ be positive integers such that $p=b^{2}+bc+c^{2}$, let $m\geq 2$ be an integer, and let $\mathcal Q$ be a finite directly regular $4$-polytope of type $\{\{6,3\}_{(m,0)},\{3,3\}\}$ such that $p(p^{2}-1)$ does not divide $|\Gamma(\mathcal{Q})|$. Then there exists a chiral $4$-polytope of type $\{\{6,3\}_{(m,0)},\{3,3\}\}$ if $p\!\mid\!m$ or $\{\{6,3\}_{(mb,mc)},\{3,3\}\}$ if $p\!\nmid\! m$, whose group is $L_{2}(p)\times \Gamma^{+}({\cal Q})$. The facets themselves are regular or chiral, respectively, in the two cases.
\end{THM}

\noindent\textbf{Proof}.
For $\mathcal P$ we now take the chiral $4$-polytope $\mathcal P$ of type $\{\{6,3\}_{(b,c)},\{3,3\}\}$ with group $L_2(p)$ described in \cite[Cor. 9.10]{SW2}. Then $\mathcal{P}\mix \mathcal{Q}$ is a chiral 
$4$-polytope with vertex-figures $\{3,3\}$ and group $L_{2}(p)\times\Gamma^{+}({\cal Q})$. Its facets are isomorphic to $\{6,3\}_{(b,c)} \mix \{6,3\}_{(m,0)}$. Now bear in mind that $[6,3]^{+}_{(b,c)} \cong 
C_p\rtimes C_6$ and $[6,3]^{+}_{(m,0)} \cong (C_m\times C_m)\rtimes C_6$, with the second factor $C_6$ in the semi-direct product coming from the rotation about the center of the base face. As in the 
previous proof, if $p\!\mid\!m$, then $\{6,3\}_{(m,0)}$ covers $\{6,3\}_{(b,c)}$ and hence the facets are isomorphic to $\{6,3\}_{(m,0)}$. On the other hand, if $p\!\nmid\!m$, then 
\[ [6,3]^{+}_{(b,c)} \mix [6,3]^{+}_{(m,0)} = [6,3]^{+}_{(mb,mc)}\cong 
(C_p \times C_m\times C_m)\rtimes C_6 , \]
of order $6pm^2$, so the facets are maps $\{6,3\}_{(mb,mc)}$.
\hfill $\square$

\begin{THM}
\label{loctor634}
Let $p$ be a prime with $p \equiv 1 \bmod 24$, let $b$ and $c$ be positive integers such that $p=b^{2}+bc+c^{2}$, let $m\geq 2$ be an integer, and let $\mathcal Q$ be a finite directly regular $4$-polytope of type $\{\{6,3\}_{(m,0)},\{3,4\}\}$ such that $p(p^{2}-1)$ does not divide $|\Gamma(\mathcal{Q})|$. Then there exists a chiral $4$-polytope of type $\{\{6,3\}_{(m,0)},\{3,4\}\}$ if $p\!\mid\!m$ or $\{\{6,3\}_{(mb,mc)},\{3,3\}\}$ if $p\!\nmid\! m$, whose group is $L_{2}(p)\times \Gamma^{+}({\cal Q})$. The facets themselves are regular or chiral, respectively, in the two cases.
\end{THM}

\noindent\textbf{Proof}.
From \cite[Thm. 2]{nosc} we obtain a chiral $4$-polytope $\mathcal P$ of type $\{\{6,3\}_{(b,c)},\{3,4\}\}$ with group $L_2(p)$. The analysis of the previous case carries over with little change and proves that the corresponding mix $\mathcal{P}\mix \mathcal{Q}$ is a chiral $4$-polytope with group $L_{2}(p)\times 
\Gamma^{+}({\cal Q})$ of the desired type.
\hfill $\square$

\begin{THM}
\label{loctor635}
Let $p$ be a prime with $p \equiv 1, 49 \bmod 60$, let $b$ and $c$ be positive integers such that $p=b^{2}+bc+c^{2}$, let $m\geq 2$ be an integer, and let $\mathcal Q$ be a finite directly regular $4$-polytope of type $\{\{6,3\}_{(m,0)},\{3,5\}\}$ such that $p(p^{2}-1)$ does not divide $|\Gamma(\mathcal{Q})|$. Then there exists a chiral $4$-polytope of type $\{\{6,3\}_{(m,0)},\{3,5\}\}$ if $p\!\mid\!m$ or $\{\{6,3\}_{(mb,mc)},\{3,5\}\}$ if $p\!\nmid\! m$, whose group is $L_{2}(p)\times \Gamma^{+}({\cal Q})$. The facets themselves are regular or chiral, respectively, in the two cases.
\end{THM}

\noindent\textbf{Proof}.
Now \cite[Thm. 6]{nosc} yields a chiral $4$-polytope $\mathcal P$ of type $\{\{6,3\}_{(b,c)},\{3,5\}\}$ with group $L_2(p)$. Then the mix $\mathcal{P}\mix \mathcal{Q}$ is again a chiral polytope with group $L_{2}(p)\times \Gamma^{+}({\cal Q})$ and its type is as stated.
\hfill $\square$

To illustrate the last three theorems, if $\mathcal{Q}$ is the universal $4$-polytope $\{\{6,3\}_{(2,0)},\{3,r\}\}$ with full group $[3,3,r] \times C_2$ (see \cite[Thm. 11B5]{arp}), then we obtain chiral $4$-polytopes of type $\{\{6,3\}_{(2b,2c)},\{3,r\}\}$, for $r=3,4,5$, with groups $L_{2}(p)\times [3,3,r]^+$ (and $p$ depending on $r$). For example, if $r=3$, the resulting polytopes have groups $L_{2}(p)\times A_{5}$.

There are variants of our constructions employing regular polytopes $\mathcal Q$ with toroidal facets $\{4,4\}_{(m,m)}$ or $\{6,3\}_{(m,m)}$, or, in some cases, other congruence classes of primes; for example, by using the polytopes of \cite[Thm. 7]{nosc} as chiral component of the mix, we can construct chiral polytopes of type $\{6,3,5\}$ with groups $L_{2}(p^2)\times \Gamma^{+}({\cal Q})$ for primes $p \equiv 7,13 \bmod 15$ such that  $p^{2}(p^{4}-1)$ does not divide $|\Gamma(\mathcal{Q})|$.

\section{Chiral 4-polytopes with toroidal facets and vertex-figures}
\label{ch4tofa}

Next we investigate locally toroidal $4$-polytopes of the remaining types $\{4,4,4\}$, $\{6,3,6\}$ and $\{3,6,3\}$. These do not have spherical vertex-figures or facets. Here we must establish the intersection property directly and then apply Lemmas~\ref{polmix} and \ref{chircritone}. We focus on self-dual polytopes, although there are similar results for general polytopes.

\begin{THM}
\label{loctor444}
Let $p$ be a prime with $p \equiv 1 \bmod 8$, let $b$ and $c$ be positive integers such that $p=b^{2}+c^{2}$, let $m$ be an odd integer, and let $\mathcal Q$ be a finite self-dual directly regular $4$-polytope of type $\{\{4,4\}_{(m,0)},\{4,4\}_{(m,0)}\}$ such that $p(p^{2}-1)$ does not divide $|\Gamma(\mathcal{Q})|$. Then there exists a properly self-dual chiral $4$-polytope of type $\{\{4,4\}_{(m,0)},\{4,4\}_{(m,0)}\}$ if $p\!\mid\!m$ or
$\{\{4,4\}_{(mb,mc)},\{4,4\}_{(mc,mb)}\}$ if $p\!\nmid\! m$, whose group is $L_{2}(p)\times \Gamma^{+}({\cal Q})$.
\end{THM}

\noindent\textbf{Proof}.
Now \cite[Cor. 8.5]{SW2} provides a properly self-dual chiral $4$-polytope $\mathcal P$ of type 
$\{\{4,4\}_{(b,c)},\{4,4\}_{(c,b)}\}$ with group $L_2(p)$. (Note that \cite[Cor. 8.5]{SW2} only states that $\mathcal P$ is self-dual but inspection of the proof shows that $\mathcal P$ is actually properly self-dual.) Mixing with $\mathcal Q$ and applying Lemmas~\ref{polmix} and \ref{chircritone} then yields a chiral $4$-polytope, $\mathcal{P}\mix \mathcal{Q}$, with the desired group, provided the intersection property holds. As before, we let $\tau_{i}=(\sigma_{i},\sigma_{i}')$ for $i=1,2,3$ be the generators of
\[ \Gamma :=\Gamma^{+}({\cal P}) \mix\Gamma^{+}({\cal Q}) = 
L_{2}(p)\times \Gamma^{+}({\cal Q}), \]
where $\Gamma^{+}({\cal P})=\langle\sigma_1,\sigma_2,\sigma_3\rangle$ and
$\Gamma^{+}({\cal Q})=\langle\sigma_{1}',\sigma_{2}',\sigma_{3}'\rangle$. Then the subgroups 
$\Gamma_{3} := \langle\tau_{1},\tau_{2}\rangle$ and 
$\Gamma_{1} := \langle\tau_{2},\tau_{3}\rangle$ for the facet or vertex-figure, respectively, of the mix $\mathcal{P}\mix \mathcal{Q}$ are themselves mixes of the facet subgroups or vertex-figure subgroups of the component polytopes. Hence it follows as before that the facets and vertex-figures themselves are either both isomorphic to $\{4,4\}_{(m,0)}$ if $p\!\mid\!m$, or isomorphic to $\{4,4\}_{(mb,mc)}$ or $\{4,4\}_{(mc,mb)}$, respectively, if $p\!\nmid\! m$.  It remains to verify the intersection property and establish self-duality.

Since the facet and vertex-figure subgroups already have the intersection property, it suffices to prove $\Gamma_{1}\cap\Gamma_{3} = \langle\tau_2\rangle$ (see Lemma~\ref{fewint}). First observe that, since $m$ and $p$ are odd and $\Gamma_{1}$ has order $4m^2$ or $4pm^2$, the order of $\Gamma_{1}\cap\Gamma_{3}$ cannot be divisible by $8$; on the other hand, $\Gamma_{1}\cap\Gamma_{3}$ does contain the subgroup $\langle\tau_2\rangle\cong C_4$.
Now, appealing to the intersection property of the component groups, if $(\varphi,\psi)\in\Gamma_{1}\cap\Gamma_{3}$, then 
\[ \varphi\in\langle\sigma_{1},\sigma_{2}\rangle \cap \langle\sigma_{2},\sigma_{3}\rangle=
\langle\sigma_{2}\rangle,\quad
\psi\in\langle\sigma_{1}',\sigma_{2}'\rangle \cap \langle\sigma_{2}',\sigma_{3}'\rangle=
\langle\sigma_{2}'\rangle,\]
and therefore $(\varphi,\psi)\in\langle\sigma_{2}\rangle\times\langle\sigma_{2}'\rangle \cong C_{4}^2$. Hence 
$\Gamma_{1}\cap\Gamma_{3}$ lies in a group isomorphic to $C_4$; as its order cannot be divisible by $8$, this forces $\Gamma_{1}\cap\Gamma_{3}= \langle\tau_2\rangle$. Thus $\Gamma$ must have the intersection property.

Finally, since $\mathcal P$ is properly self-dual, both groups $\Gamma^{+}({\mathcal P})$ and $\Gamma^{+}({\mathcal Q})$ admit involutory group automorphisms mapping $\sigma_j$ to $\sigma_{n-j}^{-1}$ and $\sigma_j'$ to $(\sigma_{n-j}')^{-1}$ for each $j$, respectively (see Section~\ref{polgr}). Hence their mix admits an involutory group automorphism mapping $\tau_j$ to $\tau_{n-j}^{-1}$ for each $j$. Thus $\mathcal{P}\mix \mathcal{Q}$ is also properly self-dual.
\hfill $\square$

In Theorem~\ref{loctor444}, if $\mathcal{Q}$ is the universal $4$-polytope $\{\{4,4\}_{(3,0)},\{4,4\}_{(3,0)}\}$ with full group $S_{6}\times C_{2}$ (see \cite[p.371]{arp}, then we arrive at properly self-dual chiral $4$-polytopes of types $\{\{4,4\}_{(3b,3c)},\{4,4\}_{(3c,3b)}\}$ with groups $L_{2}(p)\times A_{6}\times C_{2}$.

Note that the toroidal vertex-figure of the chiral polytope of Theorem~\ref{loctor444} is the enantiomorphic copy of the facet (if $p\!\nmid\! m$); that is, the two coordinates of the subscript vectors are switched. The latter will not remain true for the next two Schl\"afli types.

\begin{THM}
\label{loctor363}
Let $p$ be a prime with $p \equiv 1 \bmod 12$, let $b$ and $c$ be positive integers such that $p=b^{2}+bc+c^{2}$, let $m$ be an odd integer with $3\!\nmid\! m$, and let $\mathcal Q$ be a finite self-dual directly regular $4$-polytope of type $\{\{3,6\}_{(m,0)},\{6,3\}_{(m,0)}\}$ such that $p(p^{2}-1)$ does not divide $|\Gamma(\mathcal{Q})|$. Then there exists a properly self-dual chiral $4$-polytope of type $\{\{3,6\}_{(m,0)},\{6,3\}_{(m,0)}\}$ if $p\!\mid\!m$ or $\{\{3,6\}_{(mb,mc)},\{6,3\}_{(mb,mc)}\}$ if $p\!\nmid\! m$, whose group is $L_{2}(p)\times \Gamma^{+}({\cal Q})$.
\end{THM}

\noindent\textbf{Proof}.
We proceed exactly as in the proof of the previous theorem. Here \cite[Cor. 10.5]{SW2} provides the initial properly self-dual chiral polytope $\mathcal P$ of type $\{\{3,6\}_{(b,c)},\{6,3\}_{(b,c)}\}$ with group $L_{2}(p)$, which then is mixed with $\mathcal Q$. The corresponding group $L_{2}(p)\times \Gamma^{+}({\cal Q})$ again has the intersection property, for reasons very similar to those in the previous proof. Now $\Gamma_{1}\cap\Gamma_{3}$ lies in a group isomorphic to $C_6^2$; on the other hand, by our assumptions on $p$ and $m$, its order cannot be divisible by $4$ or $9$, so it must be $6$.
\hfill $\square$

Relatively little is known about regular $4$-polytopes of type $\{3,6,3\}$ (see \cite[Ch.11E]{arp}); in particular, the universal polytopes $\{\{3,6\}_{(m,0)},\{6,3\}_{(m,0)}\}$ have not yet been fully classified (but it is known that the polytope is finite if $m=3$). On the other hand, for every prime $m>3$ there are regular $4$-polytopes of type $\{\{3,6\}_{(m,0)},\{6,3\}_{(m,0)}\}$ with full automorphism groups given by certain finite orthogonal groups (see \cite[p.345]{moschtwo}). When these polytopes are chosen as $\mathcal{Q}$ in Theorem~\ref{loctor363}, we obtain chiral $4$-polytopes of types $\{\{3,6\}_{(mb,mc)},\{6,3\}_{(mb,mc)}\}$ as long as $p$ does not divide $m^{2}\pm 1$.

\begin{THM}
\label{loctor636}
Let $p$ be a prime with $p \equiv 1 \bmod 12$, let $b$ and $c$ be positive integers such that $p=b^{2}+bc+c^{2}$, let $m\geq 2$ be an integer with $3\!\nmid\! m$, and let $\mathcal Q$ be a finite self-dual regular $4$-polytope of type $\{\{6,3\}_{(m,0)},\{3,6\}_{(m,0)}\}$ such that $p(p^{2}-1)$ does not divide $|\Gamma(\mathcal{Q})|$. Then there exists a properly self-dual chiral $4$-polytope of type $\{\{6,3\}_{(m,0)},\{3,6\}_{(m,0)}\}$ if $p\!\mid\!m$ or $\{\{6,3\}_{(mb,mc)},\{3,6\}_{(mb,mc)}\}$ if $p\!\nmid\! m$, whose group is $L_{2}(p)\times \Gamma^{+}({\cal Q})$.
\end{THM}

\noindent\textbf{Proof}.
Now the initial properly self-dual chiral polytope $\mathcal P$ of type $\{\{6,3\}_{(b,c)},\{3,6\}_{(b,c)}\}$ with group $L_{2}(p)$ is obtained from \cite[Cor. 11.5]{SW2}. The intersection property can be settled in the same manner as in the two previous theorems. Here, $\Gamma_{1}\cap\Gamma_{3}$ lies in a group isomorphic to $C_3^2$; on the other hand, by our assumptions on $p$ and $m$, its order cannot be~$9$, so it must be $3$.
\hfill $\square$

If $\mathcal{Q}$ is the universal regular $4$-polytope $\{\{6,3\}_{(2,0)},\{3,6\}_{(2,0)}\}$ with full group $S_5\times C_2^2$ (see \cite[p.412]{arp}), then Theorem~\ref{loctor636} provides chiral $4$-polytopes of type $\{\{6,3\}_{(2b,2c)},\{3,6\}_{(2b,2c)}\}$ with groups $L_{2}(p)\times A_{5}\times C_{2}^2$.

\section{Chiral polytopes of rank $5$}
\label{ch5}

Very little is known about chiral polytopes of ranks larger than $4$. In particular, finite examples are known only for rank $5$ and were discovered quite recently in \cite{chp} (see also our concluding remarks).  Infinite examples had been found earlier in \cite{SW3} by establishing the following extension theorem and applying it to suitable polytopes of rank $4$.  {\em If $\mathcal K$ is a chiral $n$-polytope with regular facets $\mathcal F$, then there exists a chiral $(n+1)$-polytope $\mathcal L$ with facets isomorphic to $\mathcal K$; moreover, among all chiral $(n+1)$-polytopes with facets isomorphic to $\mathcal K$, there exists a universal such polytope, whose group is a certain amalgamated product of $\Gamma(\mathcal{K})$ and $\Gamma(\mathcal{F})$, with amalgamation along two subgroups isomorphic to $\Gamma^{+}(\mathcal{F})$\/}.  Note here that the assumption on the regularity on $\mathcal{F}$ is necessary, since the $(n-1)$-faces of a chiral $(n+1)$-polytope must be regular (see \cite{SW1}).

In this section we employ our methods to construct more chiral polytopes for rank $5$. Now, unlike in previous sections, the two components of the mix will not have the same Schl\"afli symbol. In particular, we establish the following theorem which enables us to construct many new chiral $n$-polytopes from a given chiral $n$-polytope, for any rank $n\geq 3$. 

\begin{THM}
\label{rankfree}
Let $\mathcal{P}$ be a finite chiral $n$-polytope of type $\{p_{1},\ldots,p_{n-1}\}$, let $\mathcal{Q}$ be a finite directly regular $n$-polytope of type $\{q_{1},\ldots,q_{n-1}\}$, and let $p_{j},q_{j}$ be relatively prime integers for each $j=1,\ldots,n-1$. Then $\mathcal{P}\mix \mathcal{Q}$ is a chiral or directly regular $n$-polytope of type $\{p_{1}q_{1},\ldots,p_{n-1}q_{n-1}\}$ with group $\Gamma(\mathcal{P})\times\Gamma^{+}(\mathcal{Q})$. Moreover, if $\mathcal{P}$ is totally chiral and 
$\Gamma(\mathcal{P})\times\Gamma^{+}(\mathcal{Q})$ does not have a subgroup which has a quotient
$\Gamma(\mathcal{P})\times\Gamma^{+}(\mathcal{P})$ (in particular if $|\Gamma(\mathcal{P})|$ does not divide $|\Gamma^{+}(\mathcal{Q})|$), then $\mathcal{P}\mix \mathcal{Q}$ is chiral.
\end{THM}

\noindent\textbf{Proof}.
Consider the mix $\Gamma:=\Gamma(\mathcal{P})\mix \Gamma^{+}(\mathcal{Q})$ of the group of $\mathcal{P}$ and the rotation subgroup of $\mathcal{Q}$. Let $\tau_{j}=(\sigma_{j},\sigma_{j}')$ for $j=1,\ldots, n-1$ denote the distinguished generators of $\Gamma$, where again 
$\Gamma(\mathcal{P}) = \langle\sigma_{1},\ldots,\sigma_{n-1}\rangle = \Gamma^{+}(\mathcal{P})$ and $\Gamma^{+}(\mathcal{Q}) = \langle\sigma_{1}',\ldots,\sigma_{n-1}'\rangle$. Since $p_{j},q_{j}$ are relatively prime, we have
\[ \langle\tau_{j}\rangle = \langle \sigma_{j}\rangle \times \langle \sigma_{j}'\rangle
= C_{p_{j}} \times C_{q_{j}} = C_{p_{j}q_{j}},\]
so in particular, $(\sigma_{j},\epsilon), (\epsilon,\sigma_{j}') \in \Gamma$ for each $j$. Then it follows immediately that, for any $J\subset \{1,\ldots,n-1\}$, 
\begin{equation}
\label{taus} 
\langle \tau_{j} \mid j\in J\rangle =  
\langle \sigma_{j} \mid j\in J\rangle \times \langle \sigma_{j}' \mid j\in J\rangle .
\end{equation}
In particular, $\Gamma = \Gamma(\mathcal{P})\times\Gamma^{+}(\mathcal{Q})$.

Now the proof of the intersection property of $\Gamma$ is quite simple and applies Lemma~\ref{fewint} inductively as follows. Define $\Gamma_{\leq k}:=\langle \tau_{j} \mid j\leq k\rangle$ for $k=1,\ldots,n-1$, so  in particular $\Gamma_{\leq n-1}=\Gamma$. Then clearly $\Gamma_{\leq 1}$ has the intersection property.

Now suppose that $\Gamma_{\leq k}$ has the intersection property for some $k\leq n-2$. We claim that then $\Gamma_{\leq k+1}$ must also have the intersection property. In fact, we have
\[ \begin{array}{rcl}
\Gamma_{\leq k} \,\cap\, \langle\tau_{2},\ldots,\tau_{k+1}\rangle& = &
\langle\sigma_{1},\ldots,\sigma_{k}\rangle \times \langle\sigma_{1}',\ldots,\sigma_{k}'\rangle 
\;\cap\;
\langle\sigma_{2},\ldots,\sigma_{k+1}\rangle \times \langle\sigma_{2}',\ldots,\sigma_{k+1}'\rangle \\[.03in] 
&=& \langle\sigma_{2},\ldots,\sigma_{k}\rangle \times \langle\sigma_{2}',\ldots,\sigma_{k}'\rangle \\[.03in]
&=& \langle\tau_{2},\ldots,\tau_{k+1}\rangle ,
\end{array} \]
where the first and last equalities are special cases of (\ref{taus}), and the second equality follows directly from the intersection property of the groups $\Gamma(\mathcal{P})$ and $\Gamma^{+}(\mathcal{Q})$, bearing in mind that $(\sigma_{j},\epsilon), (\epsilon,\sigma_{j}') \in \Gamma$ for each $j$.
Hence, by Lemma~\ref{fewint}, $\Gamma_{\leq k+1}$ also satisfies the intersection property. 

It follows that $\Gamma=\Gamma_{\leq n-1}$ has the intersection property. Thus the mix $\mathcal{P}\mix \mathcal{Q}$ of $\mathcal{P}$ and $\mathcal{Q}$ is a chiral or directly regular $n$-polytope of type $\{p_{1}q_{1},\ldots,p_{n-1}q_{n-1}\}$ with group $\Gamma$. 

Finally, if $\mathcal{P}$ is totally chiral and $\Gamma(\mathcal{P})\times\Gamma^{+}(\mathcal{Q})$ does not have a subgroup with quotient $\Gamma(\mathcal{P})\times\Gamma^{+}(\mathcal{P})$ (in particular if $|\Gamma(\mathcal{P})|$ does not divide $|\Gamma^{+}(\mathcal{Q})|$), then Lemma~\ref{chircritone} implies that $\mathcal{P}\mix \mathcal{Q}$ is chiral. 
\hfill $\square$

The polytopes obtained from the previous theorem generally have Schl\"afli symbols with rather large entries. We conclude this section with some interesting examples of chiral polytopes of rank~$5$, but note that the method  works more generally for any rank for which suitable chiral polytopes are known, so in particular also for ranks $3$ and $4$. 

First suppose that $\mathcal{P}$ is the chiral $5$-polytope 
\[ \{\{\{3,4\},\{4,4\}_{(2,1)}\},\{\{4,4\}_{(2,1)},\{4,3\}\}\} \]
with automorphism group $S_{6}$ described in \cite[p. 126, 127]{chp}. This is the universal $5$-polytope with chiral facets $\{\{3,4\},\{4,4\}_{(2,1)}\}$ and chiral vertex-figures $\{\{4,4\}_{(2,1)},\{4,3\}\}$, and is of type $\{3,4,4,3\}$. Its group $S_{6}$ has only one non-trivial normal subgroup, $A_6$. Hence, if $\mathcal{P}$ is not totally chiral, then its chirality group $X(\mathcal{P})$ must necessarily be isomorphic to $A_6$. Now employing the technique described at the end of Section~\ref{chgroups}, computation with GAP~\cite{gap} has shown that the chirality group is indeed $A_6$. Thus $\mathcal{P}$ is not totally chiral, and we cannot simply appeal to Theorem~\ref{rankfree} to conclude that a mix of $\mathcal{P}$ with a suitable directly regular $5$-polytope is chiral.

However, under certain circumstances we can still establish the chirality of the mix by hand, even though the chiral component is not totally chiral; more precisely we will apply the criterion of Lemma~\ref{wordcrit}. Now let $\mathcal{Q}$ be a regular cubic $5$-toroid $\{4,3,3,4\}_{(s^k,0^{4-k})}$ (of type $\{4,3,3,4\}$), with $s\geq 2$ and $k=1,2$ or $4$ (see \cite[Sect. 6D]{arp}), and let $[4,3,3,4]_{(s^k,0^{4-k})}^{+}$ denote the rotation subgroup of its full automorphism group (of index $2$). Then Theorem~\ref{rankfree} gives an infinite series of chiral or directly regular $5$-polytopes $\mathcal{P}\mix \mathcal{Q}$ of type $\{12,12,12,12\}$, whose automorphism groups are of the form $[4,3,3,4]_{(s^k,0^{4-k})}^{+} \times S_6$, of orders $138\,240s^4$, $276\,480s^4$ or $1\,105\,920s^4$ as $k=1,2$ or~$4$. We need to reject the possibility that $\mathcal{P}\mix \mathcal{Q}$ is directly regular. Here we can exploit the fact that a word in the generators, and its enantiomorphic image, must represent elements of the same period in the group (see Lemma~\ref{wordcrit}).

As before, let 
$\Gamma(\mathcal{P})=\langle\sigma_{1},\ldots,\sigma_{n-1}\rangle$, 
let 
$\Gamma^{+}(\mathcal{Q}) = \langle\sigma_{1}',\ldots,\sigma_{n-1}'\rangle$, 
and let 
$\Gamma := \Gamma(\mathcal{P})\times\Gamma^{+}(\mathcal{Q}) = \langle\tau_{0},\ldots,\tau_{n-1}\rangle$, with $\tau_{j}=(\sigma_{j},\sigma_{j}')$ for each $j$. Recall that $\Gamma(\mathcal{P})$ is the quotient of the Coxeter group $[3,4,4,3]$ defined by the single extra relation 
$(\sigma_{2}^{-1} \sigma_{3})^{2}\sigma_{2}\sigma_{3}^{-1} = \epsilon$. 
This extra relation determines $[4,4]_{(2,1)}$ as a quotient of $[4,4]$ (see \cite{cm}). Now consider the corresponding word 
\[ \omega:= (\tau_{2}^{-1} \tau_{3})^{2}\tau_{2}\tau_{3}^{-1} = (\omega_{1},\omega_{2}) \] 
in $\Gamma$, whose component words are 
$\omega_{1}=(\sigma_{2}^{-1} \sigma_{3})^{2}\sigma_{2}\sigma_{3}^{-1}$ 
in $\Gamma(\mathcal{P})$ (representing the trivial element) and 
$\omega_{2}=({\sigma_{2}'}^{-1} \sigma_{3}')^{2}\sigma_{2}'{\sigma_{3}'}^{-1}$ 
in $\Gamma^{+}(\mathcal{Q})$. In particular, $\omega_{2}=\sigma_{2}'$ in $\Gamma^{+}(\mathcal{Q})$. 
In fact, since $\omega_{2}$ lies in $\langle\sigma_{2}',\sigma_{3}'\rangle \cong [3,3]^{+}$ (bear in mind that $\mathcal{Q}$ is of type $\{4,3,3,4\}$), the elements $\sigma_{2}'$ and  ${\sigma_{2}'}^{-1} \sigma_{3}'$ have  period $3$ and so 
\[ \omega_{2} = ({\sigma_{2}'}^{-1} \sigma_{3}')^{-1}\sigma_{2}'{\sigma_{3}'}^{-1}
= {\sigma_{3}'}^{-1} \sigma_{2}'\sigma_{2}'{\sigma_{3}'}^{-1} 
= {\sigma_{3}'}^{-1} {\sigma_{2}'}^{-1}{\sigma_{3}'}^{-1} 
= (\sigma_{2}' {\sigma_{3}'})^{-1}{\sigma_{3}'}^{-1} 
= (\sigma_{2}' {\sigma_{3}'}){\sigma_{3}'}^{-1} = \sigma_{2}'.\]
Thus $\omega$ has period $3$ in $\Gamma$. Now consider the enantiomorphic word of $\omega$ in $\Gamma$ given by 
\[ \overline{\omega} = (\tau_{2}^{-1}\tau_{1}^{-2}\tau_{3})^{2}\tau_{1}^{2}\tau_{2}\tau_{3}^{-1} 
= (\overline{\omega_{1}},\overline{\omega_{2}}), \]
where $\overline{\omega_{1}}$ and $\overline{\omega_{2}}$ are the enantiomorphic words of $\omega_1$ and $\omega_2$ in their groups, respectively. Then a simple computation with GAP (or by hand, using the permutation representation of $\Gamma(\mathcal{P})$ in \cite[p. 126]{chp}) shows that $\overline{\omega_{1}}$ has period $5$ in $\Gamma(\mathcal{P})$. Hence $\overline{\omega}$ cannot have period $3$ in $\Gamma$. Thus the mix $\mathcal{P}\mix \mathcal{Q}$ must be chiral, since $\omega$ and $\overline{\omega}$ have different periods.

However, we can also obtain chiral polytopes of rank $5$ by directly appealing to Theorem~\ref{rankfree}. In fact, with the same choice of $\mathcal{Q}$, the two chiral $5$-polytopes of \cite[p. 127, 128]{chp} of type $\{3,8,8,3\}$ and with groups $A_{12}$ or $A_{20}$, respectively, give two infinite series of chiral $5$-polytopes of type $\{12,24,24,12\}$, with groups $A_{12}\times [4,3,3,4]_{(s^k,0^{4-k})}^{+}$ and $A_{20}\times [4,3,3,4]_{(s^k,0^{4-k})}^{+}$. Now the two initial chiral $5$-polytopes of type $\{3,8,8,3\}$ are totally chiral since their groups are simple.

We can also establish some interesting results about chiral polytopes (of ranks $3$, $4$ or $5$) with preassigned Schl\"afli symbols. For example, let $k,l,m,n \geq 3$, with $k,n$ each relatively prime to $3$ and with $l,m$ odd. Suppose that there are infinitely many directly regular $5$-polytopes of type $\{k,l,m,n\}$. Then there are also infinitely many chiral $5$-polytopes of type $\{3k,8l,8m,3n\}$. The latter are obtained  directly from Theorem~\ref{rankfree} by mixing the regular $5$-polytopes with the above (totally chiral) chiral $5$-polytopes of type $\{3,8,8,3\}$. 
\bigskip

\begin{center}
\bf Concluding remarks and acknowledgment
\end{center}
The present paper, by itself, does not establish the existence of chiral polytopes of any rank. In fact, the method already requires as input the existence of certain chiral polytopes of the specified rank. 

Most of our results were obtained by January 2009. Since then, two important developments have occurred independent of our work. We have heard from Marston Conder and Alice Devillers that they succeeded in constructing chiral polytopes of ranks $6$, $7$ and $8$ (see \cite{cd}), and Daniel Pellicer has announced that he can construct chiral polytopes of any rank $n\geq 3$ (see \cite{pel}). A new example of rank $5$ is also described in Pellicer-Weiss~\cite{pelwei}. There is a good chance that some of these polytopes can serve as input polytopes in our construction and hence provide many more chiral polytopes of the specified ranks.

This research was begun while the second and third authors visited the first at the University of Aveiro in April 2008. We would like to thank the University of Aveiro for its hospitality.

\end{document}